\documentclass[hidelinks,onefignum,onetabnum]{siamart251216}

\usepackage{lipsum}
\usepackage{amsfonts}
\usepackage{graphicx}
\usepackage{epstopdf}
\usepackage{algorithmic}
\usepackage{amsopn}
\usepackage{enumerate,bm}
\usepackage{hyperref}
\usepackage{xcolor}
\usepackage{soul}
\usepackage{multirow}
\usepackage{amssymb,amsmath}
\usepackage{tensor}
\usepackage{cancel}
\usepackage{mdframed}
\usepackage{makecell}
\usepackage{booktabs}
\usepackage{subcaption}
\usepackage{graphbox}
\usepackage{float}

\newsiamremark{remark}{Remark}
\newsiamremark{hypothesis}{Hypothesis}
\crefname{hypothesis}{Hypothesis}{Hypotheses}
\newsiamthm{claim}{Claim}
\newsiamremark{fact}{Fact}
\crefname{fact}{Fact}{Facts}

\headers{Optical Fourier Transform-Based Preconditioning}{}

\title{Optical Fourier Transform-Based Preconditioning for Large-Scale Inverse Problems}

\author{C. Drum\thanks{Department of Mathematics, Emory University, Atlanta, GA
    (\email{cndrum@emory.edu},\\ \email{jnagy@emory.edu}).}
  \and L. Onisk\footnotemark[1]\,\,\,\thanks{Argonne Leadership Computing Facility, Argonne National Laboratory, Lemont, IL (\email{lonisk@anl.gov}).}
  \and A. Brown\footnotemark[1]\,\,\,\thanks{Mathematics Department, Spelman College, Atlanta, GA (\email{arianabrown@spelman.edu}).}
  \and J. G. Nagy\footnotemark[1]}

\begin{document}

\maketitle

\begin{abstract}
Fast Fourier Transform (FFT)-based preconditioners provide an efficient approach for approximating inverse operators in large-scale linear inverse problems when the underlying operator structure can be exploited in the Fourier domain. Optical Fourier transforms (OFTs), the analog of digital Fourier transforms, exploit light propagation through optical systems and offer the potential for greater throughput and energy efficiency than conventional digital FFT implementations. However, practical optical systems are inherently imperfect and introduce a variety of errors into the computed output. Thus, OFTs can only provide approximations to their digital counterpart, and the implications of these inexact transforms for iterative solvers remains largely unexplored. This work investigates the feasibility of OFT-based preconditioning of iterative Krylov subspace methods for the solution of discrete inverse problems. We digitally simulate the optical Fourier transform using the FFT and incorporate sources of error expected in practical optical systems into the simulation. Numerical experiments suggest that given sufficient optical system quality, iterative methods utilizing OFT-based preconditioners can maintain reconstruction quality and convergence comparable to those with fully digital FFT-based preconditioners.

\end{abstract}



\section{Introduction} \label{sec:intro}

We are interested in computing approximate solutions to linear least-squares problems of the form
\begin{equation} \label{leastSquares}
\min_{x \in \mathbb{R}^n} \|Ax-b\|
\end{equation}
using efficient preconditioning strategies that are enabled by hybrid digital-optical processes. We denote the Euclidean norm throughout by $\|\cdot\|$. Here, $A \in \mathbb{R}^{m \times n}$ is a matrix whose singular values decay without significant gap and cluster at numerical zero (i.e., the matrix is ill-conditioned). The methods discussed herein allow for $m \geq n$ as well as $m<n$. Problems of this type are usually referred to as linear discrete ill-posed problems. They naturally arise in massive data streaming problems such as the training of the random feature model in machine learning \cite{RB07}, or limited angle tomography problems such as those from medical imaging \cite{Chung_inter3,streamInvProb}. They also classically arise from the discretization of Fredholm integral equations of the first kind \cite{engl1996regularization,Ha98}.

The data vector $b$ from \eqref{leastSquares} in many inverse problems often represents measurements that are corrupted by noise contamination or instrumental disturbances. Numerically, this error may also include truncation or approximation errors. We collectively denote these errors by $e\in \mathbb{R}^m$ so that
\begin{equation*}
    b = \hat{b}+e
\end{equation*}
where $\hat{b}$ is the unknown and error-free vector associated with $b$. Ideally, we would like to determine the solution of the minimal Euclidean norm, $x^{\dagger}$, associated with the unavailable least-squares problem 
\begin{equation}
    x^{\dagger} = \text{arg}\min_{x \in \mathbb{R}^n}\left\|Ax-\hat{b} \right\|.
\end{equation}
However, it is well understood that because the available data $b$ is contaminated by error and the singular values of $A$ cluster near the origin, that the solution of the available problem \eqref{leastSquares} is usually a poor approximation to $x^{\dagger}$. To derive a meaningful solution, regularization is often employed.

In variational regularization, prior knowledge of the solution and its properties (e.g., smoothness or sparsity) is included by replacing \eqref{leastSquares} with a penalized least-squares problem. Herein, we focus on standard Tikhonov regularization, 
\begin{equation}  \label{eq:3} 
\min_{x \in \mathbb{R}^n} \left\{\left\| Ax-b \right\|^{2} + \alpha^{2} \left\|x\right\|^{2}\right\}, 
\end{equation} 
where $\alpha > 0$ is a regularization parameter which may be thought of as a ``trade-off" parameter that balances the sensitivity of the computed solution to the error in $b$, as well as the closeness in norm to the desired solution $x^{\dagger}$. An equivalent form of \eqref{eq:3} is the regularized normal equations
\begin{equation} \label{eq:NormalEqns}
    \left(A^TA+\alpha^2I\right)x=A^Tb.
\end{equation}

Approximating the solution to the Tikhonov problem \eqref{eq:3} is usually untenable via direct methods due to its size and because the forward model matrix $A$ (and its adjoint) are often only accessible through matrix-vector products. Instead, we consider the mathematically equivalent saddle-point formulation given by
\begin{equation} \label{eq:saddlePoint}
\underbrace{\begin{bmatrix}
I & A \\ A^T  & \alpha^2 I 
\end{bmatrix}}_K
    \underbrace{\begin{bmatrix}
        r \\ x
    \end{bmatrix}}_{s}
    =
    \underbrace{\begin{bmatrix}
        b \\ 0
    \end{bmatrix}}_{c}
\end{equation}
where $r=b-Ax\in \mathbb{R}^m$, $K \in \mathbb{R}^{(m+n) \times (m+n)}$, $s \in \mathbb{R}^{m+n}$, and $c \in \mathbb{R}^{m+n}$. The first block row $r+Ax=b$ is the residual and the second block row $A^Tr-\alpha^2x=0$ represents a reformulation of the regularized normal equations \eqref{eq:NormalEqns}. In this form, the residual and computed solution are calculated simultaneously; however, the matrix associated with the regularized normal equations \eqref{eq:NormalEqns} is not formed explicitly. Rather than computing the direct solution to \eqref{eq:saddlePoint}, iterative strategies such as applying Krylov subspace methods may be utilized since they only require forward operations with $A$ and $A^T$ \cite{benzi2005numerical}. However, because it can take many iterations for standard methods (e.g., CGLS, GMRES, etc.) to capture the dominant basis information necessary to constitute a meaningful approximate solution for inverse problems, preconditioning is often employed to reduce the number of iterations necessary.

Discrete variants of the Fourier transform (FT) have long been utilized as preconditioners for discrete ill-posed problems; e.g., \cite{HankeNagy1996, HankeNagyPlemmons1993, Hansen2002, NagyPlemmonsTorgersen1996}. Classically, image deconvolution problems whose system matrices are represented by Toeplitz matrices have been approximated by circulant matrices that can be diagonalized by the fast Fourier transform (FFT); see e.g., \cite{chan1993fft, chan1994circulant, DonatelliHanke2013, van1993approximation}. The most widely used FFT algorithm, due to Cooley and Tukey, contains a computational complexity of $\mathcal{O}(n\log n)$ which, although efficient, suffers from constraints due to digital hardware, including memory bandwidth and energy use. For exceptionally large computational workloads, data movement and memory access can dominate the cost, often exceeding the cost of floating point operations and becoming the primary bottleneck for performance. These challenges are especially relevant in optimizing the performance of modern exascale computing environments, see, e.g., \cite{allen2025aurora,frontier2023}.

Our particular focus herein is on determining the experimental viability of utilizing optical Fourier transforms (OFTs) within preconditioners to compute an iterative approximation to \eqref{eq:saddlePoint}. In contrast to digital FFTs, which require a sequence of arithmetic operations and data movement through electronic memory hierarchies, OFTs instead operate through diffraction and light propagation, allowing the transformation to occur at the speed at which light propagates through the optical system. Furthermore, the natural parallelism of OFTs makes them appealing for high-throughput computation with the potential for reduced energy consumption. Various physical implementations of optical systems have been developed, including lens-based free-space systems and photonic integrated circuits \cite{Kazanskiy2026, McMahon_2023}. In this work, we focus on lens-based systems as they are well-established and supported by extensive literature detailing their theoretical behavior as well as the various sources of error that arise in practical systems; see, e.g., \cite{BornWolf, Goodman, Greivenkamp,Voelz}.

Optical analog computing has received recent interest in scientific computing due to the ability of various applications to tolerate limited precision and hardware-induced noise, including, for example, neural networks \cite{Chang_OpticalNN_2018,OpticalNN_21,Liu_24,Meng_2025,Spall_2022}. However, the use of optical computing within iterative numerical algorithms, particularly with respect to preconditioning to the best of our knowledge has been unexplored. Recent work utilizing hybrid digital-analog computing approaches focused on resistive crossbar arrays has shown promise for preconditioning; however, feasible array sizes remain relatively small and require block-based approaches for larger problems \cite{Feinberg2021, Shah_digitalAnalog}. We begin our investigation by providing a model of the OFT with respect to the Fourier operator that incorporates physically motivated sources of error encountered in practical optical systems. To investigate the feasibility of using OFTs within preconditioners, we simulate a hybrid digital-optical framework by including these errors in the discrete FFT and numerically evaluate the preconditioner by utilizing it within an iterative algorithm applied to \eqref{eq:saddlePoint} on two discrete ill-posed inverse test problems.

The remainder of the paper is structured as follows. In Section \ref{sec:background}, we review the iterative generalized minimal residual (GMRES) method and its flexible variant that we will use to approximate the solution to \eqref{eq:saddlePoint}. We also define the block preconditioning types that allow for natural integration with our simulated OFT. Section \ref{sec:oft} provides an overview of the continuous OFT and the errors we consider throughout that are realized in practical OFT systems. In Section \ref{sec:ImpSim}, we describe the simulation of our modeled OFT and the incorporation of practically realized quantization errors. The numerical results and relevant preliminaries are provided in Section \ref{sec:results}. Concluding remarks and an outlook are given in Section \ref{sec:conclusions}.

\section{Background} \label{sec:background}
This section provides necessary background for approximating the solution to the saddle-point problem \eqref{eq:saddlePoint} using Krylov subspace methods. Specifically, we provide a brief overview of the GMRES method and its flexible variant, FGMRES. We then consider and discuss two preconditioner types whose block structures integrate naturally with the OFT. 

\subsection{Preconditioned System and GMRES/FGMRES Methods} \label{sec:krylov}

To investigate the OFT as an effective preconditioner we begin by considering the preconditioning of the minimization problem associated with \eqref{eq:saddlePoint}. Since we assume that we only have access to forward operations with $A$ and $A^T$ and do not wish to form the corresponding block augmented matrix associated with \eqref{eq:saddlePoint} explicitly; instead, we utilize the GMRES method as an iterative solver. We also consider its flexible variant, FGMRES, for the situation when the preconditioner changes at each iteration.

We consider the right preconditioned variant of the minimization problem associated with \eqref{eq:saddlePoint} given by
\begin{equation} \label{eq:min_saddleProb}
    \min_{y \in \mathbb{R}^{m+n}}\left\|KM^{-1}y-c \right\|, \qquad s=M^{-1}y
\end{equation}
where $M\in \mathbb{R}^{(m+n) \times (m+n)}$ is a preconditioner that approximates $K$ and is such that the condition number of $KM^{-1}$ is smaller than that of $K$. The GMRES method is a subspace projection technique that allows for efficient approximation to square non-symmetric linear systems \cite{saad2003iterative}. In our notation, the $p^{th}$ step of GMRES uses the Arnoldi process to compute orthonormal basis vectors of the Krylov subspace ${\cal{K}}_{p}(KM^{-1},c)$ where an approximate solution is sought. At step $p$, the Arnoldi process produces a $(p+1) \times p$ upper Hessenberg matrix, $H_{p+1}$, and forms the following relationship: 
\begin{equation} \label{eq:ArnoldiProcess}
    KM^{-1}V_{p} = V_{p+1}H_{p+1}
\end{equation} 
where $V_{p+1}\in \mathbb{R}^{(m+n) \times (p+1)}$ contains orthonormal columns and $V_p\in \mathbb{R}^{(m+n) \times p}$ contains the first $p$ columns of $V_{p+1}$ that span  ${\cal{K}}_{p}(KM^{-1},c)$; see \cite{saad2003iterative} for more details.

When $M$ can change at every iteration, a flexible variant is often more useful for iterative solvers with respect to capturing the non-deterministic sources of error of the preconditioner; we discuss these source components in Section \ref{sec:ImpSim} for the OFT. After $p$ steps, the generalized Arnoldi process results in the following relation:
\begin{equation*}
KZ_p = V_{p+1}H_{p+1}
\end{equation*} 
where $Z_p = [M_1^{-1}v_1, M_2^{-1}v_2, \cdots, M_p^{-1}v_p]\in \mathbb{R}^{(m+n) \times p}$. Unlike the standard Arnoldi process with a fixed preconditioner \eqref{eq:ArnoldiProcess}, the range of the matrix $Z_p$ containing the solution subspace does not generally contain orthonormal columns and whose span does not generally correspond to a Krylov subspace \cite{saad1993flexible}. The FGMRES algorithm is provided in Algorithm \ref{alg:FGMRES} for completeness; however, we note that the block matrices utilized within are not formed explicitly and are only shown for clarity. Additionally, we clarify that while there are alternative flexible iterative strategies including Lanczos-type algorithms that could be utilized to avoid potentially excessive basis storage, we consider GMRES-based strategies to avoid issues related to loss of orthogonality of the basis vectors in finite precision. 

\begin{algorithm}[h]
\caption{FGMRES for \eqref{eq:saddlePoint}} \label{alg:FGMRES}
\begin{algorithmic}[1]
\REQUIRE $K\in \mathbb{R}^{(m+n) \times (m+n)}$, $M\in \mathbb{R}^{(m+n) \times (m+n)}$, $c \in \mathbb{R}^{m+n}$, $s_0 \in \mathbb{R}^{m+n}$ 
\STATE Compute $r_0 = c - K s_0 $, $\beta = \|r_0\|_2$, $v_1 =  r_0/\beta$
\FOR{$p = 1,2,\ldots$} 
  \STATE $z_p = M^{-1}_pv_p$
  \STATE $w_p = Kz_p$
\FOR{$j = 1, \ldots,p$}
  \STATE $h_{j,p} = v_j^{T}w_p$; $w_p = w_p-h_{j,p}v_j$
  \ENDFOR
  \STATE $h_{p+1,p} = \|w_p\|_2 $
  \IF{$h_{p+1,p} = 0$ break}
 \STATE $v_{p+1} = w_p/h_{p+1,p}$
 \ENDIF
 \ENDFOR
  \STATE Compute $y_p = \arg \min\| \beta e_1 - H_{p+1}y \|_2^2 $ 
  \STATE $ s_p = s_0 + Z_p y_p$

\end{algorithmic} 
\end{algorithm}

\subsection{Preconditioning Structure} \label{sec:precon_types}

In this section we provide an overview of the preconditioning structures that naturally arise when considering a hybrid digital-optical framework for the OFT that we describe in Section \ref{sec:ImpSim}. Specifically, we consider two types of preconditioners $M$ that allow for efficient approximation of $K$ using matrix-vector products with the OFT.

\noindent \textbf{(Type 1)} Here, we consider the situation when the matrix $A$ inside \eqref{eq:saddlePoint} may be approximated by another matrix $C$ that may be diagonalized by the FFT, i.e., $A\approx C$, so that
\begin{equation} \label{eq:type1precon}
    K \approx M = \begin{bmatrix}
        I & C \\ C^{\ast} & \alpha^2I
    \end{bmatrix}
\end{equation}
where $^{\ast}$ denotes the conjugate transpose when $C$ is complex.
Using the Schur complement and the spectral decomposition $C=F^{\ast}\Lambda F$, application of the inverse in GMRES or FGMRES may be represented by
\begin{equation}
    M^{-1}=\begin{bmatrix}
        F^{\ast} & 0\\0 & F^{\ast}
    \end{bmatrix}
    \begin{bmatrix}
        I-\Lambda P_1\Lambda^{\ast} & \Lambda P_1\\P_1\Lambda^{\ast} & -P_1
    \end{bmatrix}
    \begin{bmatrix}
        F & 0\\0 & F
    \end{bmatrix}
\end{equation}
where $P_1=\left(\left|\Lambda\right|^2+\alpha^2I\right)^{-1}$, which can be efficiently applied since the central block matrix contains diagonal blocks. In Section \ref{sec:opt_implementation} we describe our OFT approximation to $P_1$ when utilizing a Type 1 preconditioner.

\noindent \textbf{(Type 2)} In this preconditioning structure, we consider the approximation of the block LDL factorization of $K$ given by
\begin{equation} \label{LDL_K}
    K = \underbrace{\begin{bmatrix}
        I & 0 \\A^T & I
    \end{bmatrix}}_{L}
    \underbrace{\begin{bmatrix}
        I & 0 \\0 & -\left(A^TA+\alpha^2I\right)
    \end{bmatrix}}_{D}
    \underbrace{\begin{bmatrix}
        I & A \\0 & I
    \end{bmatrix}}_{L^T}
\end{equation}
where $L,\, D \in \mathbb{R}^{(m+n) \times (m+n)}$.
To avoid computing a linear solve with $A^TA+\alpha^2I$ from \eqref{LDL_K} in the variants of GMRES we consider herein, we utilize an approximation $P_2\approx A^TA+\alpha^2I$ which may, for example, be derived from the diagonalization by the FFT. We describe our approximant to $P_2$ in Section \ref{sec:opt_implementation} when utilizing a Type 2 preconditioner. 
Structurally, we define
\begin{equation} \label{eq:type2precon}
    K^{-1} \approx M^{-1} = \begin{bmatrix}
        I & -A \\0 & I
    \end{bmatrix}
    \begin{bmatrix}
        I & 0 \\0 & -P_2^{-1}
    \end{bmatrix}
    \begin{bmatrix}
        I & 0 \\-A^T & I
    \end{bmatrix}.
\end{equation}
While saddle-point preconditioners often only utilize the approximation to the block diagonal matrix $D$ (see, e.g., \cite{benzi2005numerical}), we found that including the block triangular factors improved the robustness of the preconditioner while also taking advantage of the efficient application of matrix-vector products with $A$ and $A^T$.

\section{Optical Fourier Transform} \label{sec:oft}

In this section we provide an overview of the OFT and its physically realized errors that we consider for our investigation. By doing so, we introduce the components necessary to model the OFT based on the discrete FT which will be used to simulate the OFT for use within a preconditioner in our numerical experiments in Section \ref{sec:results}. We remark here that for the purpose of introducing the OFT and its errors, we initially describe the transform in the continuous setting. However, we return to the discrete setting in Section \ref{sec:ImpSim} when introducing our proposed model for numerical simulation.

\begin{figure}[H]
    \centering
    \includegraphics[width=0.4\linewidth]{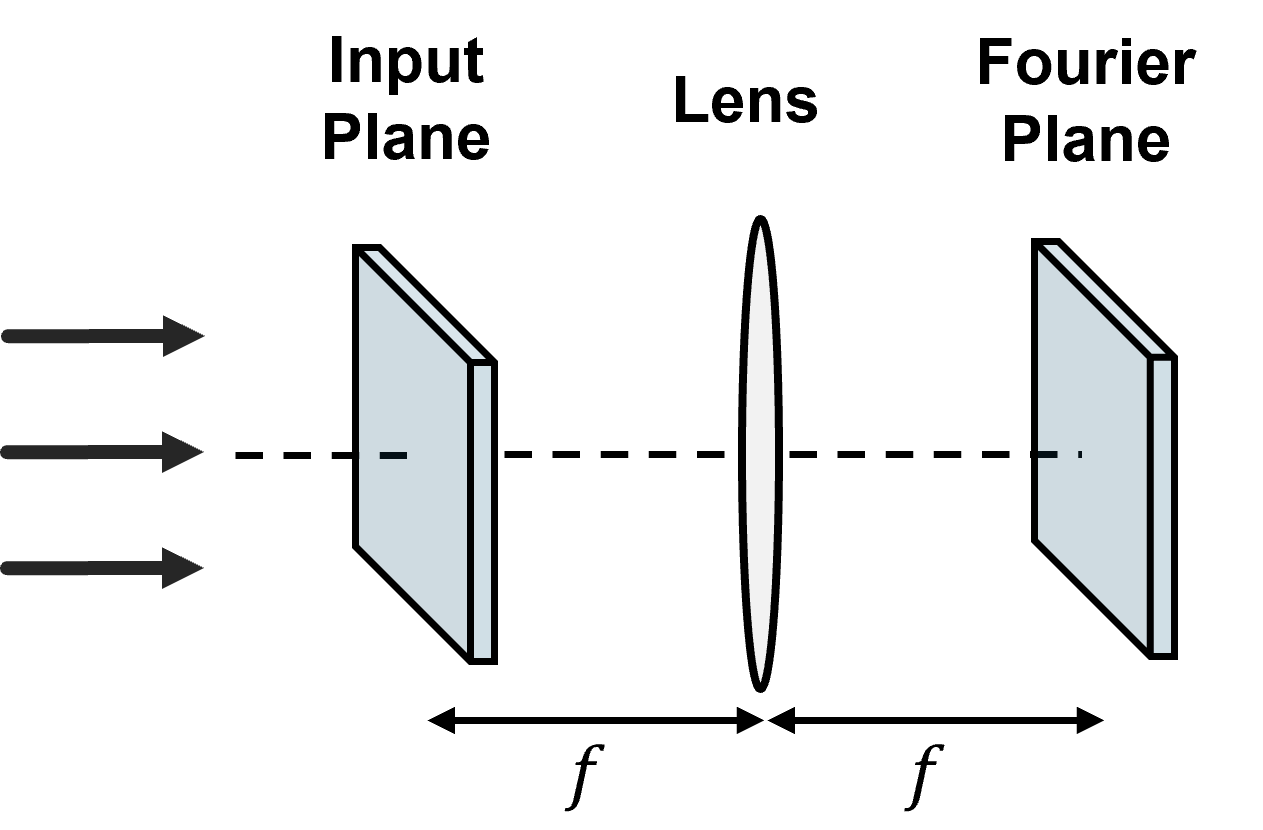}
    \caption{Schematic of a $2f$ optical system where the input field is a programmable optical light modulator that is placed a distance of one focal length $f$ in front of a thin converging lens. The Fourier transform of the input is located at the back focal plane, which is also one focal length $f$ from the lens. }
    \label{fig:2fconfig}
\end{figure}

A common configuration for the physical implementation of an OFT is a $2f$ optical system, as illustrated in Figure~\ref{fig:2fconfig}. In this configuration, the input is placed one focal length $f$ in front of a converging lens and the Fourier plane is located one focal length $f$ behind the lens. Thus, the total distance between the input and Fourier planes is $2f$. 

The input plane contains an optical modulator, such as a spatial light modulator (SLM), that modulates the optical field to encode the input data. The modulation device can be represented mathematically by the complex transmission function
$$t(\xi, \eta)=a(\xi,\eta)e^{i\phi(\xi,\eta)},$$
where $a(\xi,\eta)$ and $\phi(\xi,\eta)$ denote the amplitude transmission and phase delay, respectively. Depending on the modulation device, $t(\xi,\eta)$ may represent amplitude-only, phase-only, or complex modulation. In this work, we consider complex modulation where both the amplitude and phase of the field are modulated. Assuming coherent, monochromatic plane-wave illumination with constant complex amplitude $\mathcal{A}$, the resulting complex field immediately after the object is
\begin{equation*}
    U(\xi,\eta)=\mathcal{A}\,t(\xi,\eta).
\end{equation*}

Following Goodman \cite{Goodman}, under ideal optical conditions, when the input is placed in the front focal plane of the lens, the field in the back focal plane is given by

\begin{equation} \label{eq:OFTintegral}
U_f(u,v) = \frac{1}{i\lambda f}\iint
 U(\xi,\eta)
\exp\!\left[
-i\frac{2\pi}{\lambda f}(\xi u+\eta v)
\right]
d\xi d\eta ,
\end{equation}
where $\lambda$ denotes the wavelength, $i=\sqrt{-1}$, and $(u,v)$ are the physical coordinates in the back focal (Fourier) plane. Defining the spatial frequencies $$f_x = \frac{u}{\lambda f}, \quad \quad f_y = \frac{v}{\lambda f},$$
yields
\begin{align}
U_f(f_x,f_y)
&= \frac{1}{i\lambda f}\,\mathcal{F}\{U(\xi,\eta)\}(f_x,f_y) \nonumber\\
&= \frac{\mathcal{A}}{i\lambda f}\,\mathcal{F}\{t(\xi,\eta)\}(f_x,f_y).
\label{eqn:2f}
\end{align}
Therefore, the field in the back focal plane is proportional to the Fourier transform of the input object's transmission function. The constant factor $\frac{\mathcal{A}}{i\lambda f}$ is uniform across the Fourier plane and only changes the overall complex amplitude of the transformed field. Alternatively, when a Fourier-domain field is placed in the front focal plane, the system produces the corresponding inverse Fourier transform at the back focal plane.

\subsection{Sources of Optical Error}
Equation (\ref{eq:OFTintegral}) describes an exact Fourier transform relationship that is produced by an ideal $2f$ optical system. Practical physical implementations, however, are subject to a variety of imperfections. Here, we provide a discussion of several sources of error in optical systems, including finite aperture effects, optical aberrations, misalignment, and noise. Finite aperture effects limit the transmitted spatial frequencies, optical aberrations distort the propagated wavefront, misalignment introduces phase and positioning errors, and noise degrades the measured signal. In what follows, we describe each of these error sources and their incorporation into our model for simulating the OFT.

\subsubsection{Finite Aperture}

The aperture of an optical system is the hole or opening that limits the light being propagated through the system. All lens-based optical systems have a finite aperture due to the finite size of the lens, i.e., the lens can only transmit light within its aperture. As a result, only a limited range of spatial frequencies can be transmitted through the optical system. This leads to deviations from the ideal Fourier transform relationship. The finite extent of the lens can be accounted for by incorporating a pupil function $P(u,v)$ defined by
\[
P(u,v)=
\begin{cases}
    1 & \text{inside the lens aperture} \\
    0  & \text{otherwise}
\end{cases}
\]
and modifies the field by
\begin{equation}
    P(u,v)U_f(u,v),
\end{equation}
where $U_f(u,v)$ is the ideal Fourier plane field defined by \eqref{eqn:2f}. Thus, for a circular aperture with radius $r$, the pupil function is given by
\begin{equation}  \label{eq:pupil}
    P(u,v)=
    \begin{cases}
    1 & \sqrt{u^2+v^2}\le r \\
    0  & \sqrt{u^2+v^2}> r,
    \end{cases}
\end{equation}
where the field outside of the circle is truncated, or equivalently, the spatial frequencies outside $\displaystyle \frac{r}{\lambda f}$ are excluded from the measured Fourier transform.

Also of note is the concept of ``clear aperture'', which is defined as the region of the optical component where the manufacturer's specifications for certain optical properties (e.g., surface quality) are guaranteed to be met. Typically, the clear aperture is approximately $90\%$ or more of the full lens diameter \cite{EdmundOptics}. Outside of this region, manufacturers do not guarantee that the lens will adhere to the stated specifications.

\subsubsection{Aberrations} \label{sec:lens_abberations}

Also present in optical systems are optical aberrations, which are defined as deviations of an optical wavefront from its ideal form. These deviations can occur due to inherent properties of lenses or because of physical imperfections in the optical system. For example, focusing errors can occur due to improper positioning of optical components, while spherical aberrations are inherent to the geometry of spherical lenses. Regardless of the type or source, aberrations manifest as phase errors in the propagated optical field and can be modeled by introducing a wavefront aberration function $W(x,y)$, which represents the optical path-length error across the physical pupil of the optical system. The generalized pupil function including both finite aperture and wavefront aberration is given by
\begin{equation}
\mathcal{P}(x,y) = P(x,y)e^{ikW(x,y)},
\end{equation}
where $P(x,y)$ is the pupil function defined previously and $k=\frac{2\pi}{\lambda}$ is the wavenumber. The exponential term represents the phase distortion introduced by the aberration. Thus, aberrations are incorporated into the pupil function and introduce phase distortions into the Fourier transform. With appropriate coordinate scaling, the generalized pupil function can be expressed in spatial frequency coordinates $\mathcal{P}(f_x,f_y)$ and then the Fourier plane field in \eqref{eqn:2f} may be modified according to
\begin{equation}
    \mathcal{P}(f_x,f_y) U_f(f_x,f_y).
\end{equation}
Thus, aberrations act as multiplicative phase distortions within the region of the Fourier domain passed by the finite aperture. 

A common approach for representing wavefront aberrations is through a Zernike polynomial expansion \cite{Schmidt,Lak_Zernike}. Zernike polynomials allow the wavefront aberration function to be written as
\begin{equation} \label{eq:zernikepoly}
    W(r,\theta) = \sum_{j=1}^{m} a_j Z_j(r,\theta),
\end{equation}
where $Z_i(r,\theta)$ denotes the $i$-th Zernike polynomial and $a_i$ is the corresponding aberration coefficient. These polynomials are defined over a circular domain that represents the physical pupil of the optical system and therefore it is both convenient and common to express the coordinates in polar form. 
Although the full Zernike series is infinite, practical implementations use only a finite number of $m$ modes. These modes represent common optical aberrations, such as defocus, astigmatism, coma, and spherical aberration. For additional details on Zernike polynomials and the optical aberrations they represent; see, e.g., \cite{Mahajan_Zernike,Niu_Zernike}.

\subsubsection{Noise}

In practical optical systems, the Fourier plane signal produced by an optical Fourier transform cannot be measured exactly. The measured optical field is corrupted by noise from the physical detection process and the associated electronics. Consequently, the observed spectrum is a noisy approximation of the ideal optical field. There are many different potential sources of noise in an optical system, including detector noise, photon shot noise, and speckle noise, for instance. Each noise source could be represented using its own respective statistical noise model, however, in this work, we consider their combined effect as additive complex Gaussian noise in the Fourier plane field,
\begin{equation}
U_f(f_x,f_y)+\varepsilon(f_x,f_y),
\end{equation}
where \(\varepsilon(f_x,f_y)\) is a Gaussian random field with magnitude dependent on the conditions of the given optical system.

\subsubsection{Misalignment}

An ideal OFT assumes that the input object is located exactly in the front focal plane of the lens and that the Fourier plane detector is centered on the optical axis. In practice, small positioning errors may occur in either plane causing the resulting Fourier transform to differ from the ideal one. For our simulation purposes, we consider a lateral shift in the input plane. If the input field is displaced by $(\Delta \xi,\Delta \eta)$, then the shifted field is given by
\begin{equation}
    U(\xi-\Delta \xi,\eta-\Delta \eta),
\end{equation}
and by the Fourier shift theorem,
\begin{equation}
    \mathcal{F}\{U(\xi-\Delta \xi,\eta-\Delta \eta)\} = e^{-i2\pi(f_x\Delta \xi+f_y\Delta \eta)}\mathcal{F}\{U(\xi,\eta)\}.
\end{equation}
Thus, input plane misalignment introduces a linear phase ramp across the Fourier spectrum. Equivalently, the Fourier plane field becomes
\begin{equation} \label{eq:shift}
e^{-i2\pi(f_x\Delta \xi+f_y\Delta \eta)}
U_f(f_x,f_y),
\end{equation}
where $U_f(f_x,f_y)$ denotes the Fourier plane field corresponding to the unshifted input field $U(\xi,\eta)$. Other forms of misalignment may occur in practical optical systems, such as tilt of an optical component and axial spacing errors. However, these forms of misalignment are commonly modeled as wavefront aberrations, as described in Section~\ref{sec:lens_abberations}; see \cite{Goodman,Voelz}.

\subsection{Modeling the Optical Fourier Transform}

We now describe how the sources of optical error introduced in the preceding section are incorporated into our model for simulating the OFT. Combining the error sources yields the following model for a practical optical Fourier transform,
\begin{equation}
    \mathcal{H}(f_x,f_y)\odot U_f(f_x, f_y)+\varepsilon(f_x,f_y),
\end{equation}
where $\odot$ denotes elementwise multiplication and $\mathcal{H}$ denotes a transfer function incorporating aperture, abberation, and lateral input shift, and is defined by 
\begin{equation}
\begin{aligned}
\mathcal{H}(f_x,f_y)
&=\mathcal{P}(f_x,f_y)e^{-i2\pi(f_x\Delta \xi+f_y\Delta \eta)}\\
&=P(f_x,f_y)e^{i\left[kW(f_x,f_y)-2\pi\bigl(f_x\Delta\xi+f_y\Delta\eta\bigr)\right]}.
\end{aligned}
\end{equation}

To investigate the feasibility of OFT-based preconditioning, we simulate the OFT in the discrete setting by modeling the dominant physically motivated components of an imperfect optical system and incorporating them into a digital FFT. Specifically, let $T$ denote the discrete matrix input to the optical system. We compute the discrete Fourier transform $\mathrm{F}\{T\},$ and then apply the transfer function and additive noise according to
\begin{equation} \label{eq:ofterrs}
    \mathcal{H}\odot F\{T\}+\varepsilon
\end{equation}
Here, $\mathcal{H}$ denotes the discretized form of $\mathcal{H}(f_x,f_y)$ and $\varepsilon$ is additive noise on the discrete Fourier grid. For simplicity, we omit the explicit dependence on the spatial frequency coordinates. For the remainder of this work, we consider this discrete setting.

\section{OFT-based Preconditioning and Simulation} \label{sec:ImpSim}

We now introduce our proposed framework for implementing and simulating OFT-based preconditioning, where the central idea is to use analog OFTs in place of the digital Fourier transform operations required within FFT-based preconditioners. In Section~\ref{sec:hybrid_framework}, we first introduce the hybrid digital-optical structure and describe the data flow between the digital and optical components, the required data conversions, and the resulting quantization errors. Then, in Section~\ref{sec:opt_implementation}, we discuss the specific optical implementations for the preconditioners considered in this work, as well as their complete simulation models that incorporate both optical errors and data-conversion effects.

\subsection{Hybrid Digital-Optical Framework} \label{sec:hybrid_framework}

\begin{figure}[H]
     \centering
         \includegraphics[width=.6\textwidth]{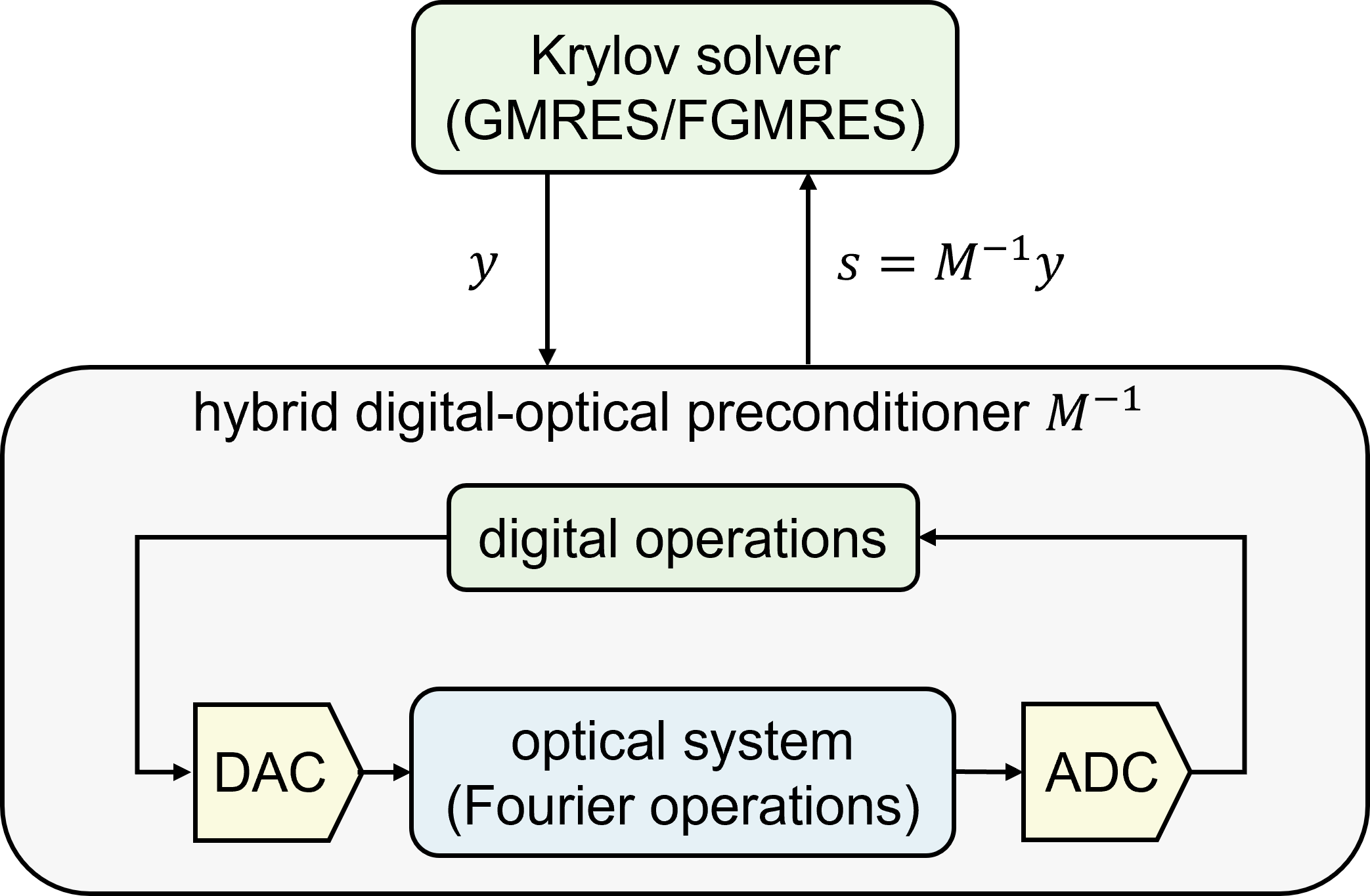}
        \caption{Data flow in the proposed hybrid digital-optical preconditioning framework. The preconditioner consists of both digital operations and optical Fourier-based operations. A DAC is required for digital-to-analog conversion between the digital and optical components. An ADC is required for analog-to-digital conversion between the optical and digital components.}
        \label{fig:hybrid_flowchart}
\end{figure}

The proposed hybrid framework offloads the Fourier operations required by a preconditioner of a linear system to an optical system, while all remaining computations are performed digitally. An important and practical consideration in hybrid digital-optical methods is how data is transferred between the digital solver and the optical system, as well as the errors that can come with this process. A digital-to-analog converter (DAC) must be used to convert the digital data to analog data for optical processing and an analog-to-digital converter (ADC) is needed to digitize the optical system output data, as illustrated in Figure~\ref{fig:hybrid_flowchart}.
Both the DAC and the ADC introduce quantization errors during the data conversion process. The DAC introduces quantization by mapping the digital values to a finite number of discrete levels that the SLM can apply. The ADC, on the other hand, must take the continuous analog signal produced by the detector at the output plane of the optical system and map it to a finite number of discrete digital values.

To further support the realism of the OFT simulation, DAC and ADC quantization are incorporated into the error model. Specifically, we utilize uniform $b$-bit quantization, described by the operator  
\begin{equation*}
Q(x) = \Delta \left\lfloor \frac{x}{\Delta} + \frac{1}{2}\right\rfloor 
\end{equation*} 
where $\displaystyle \Delta $ is the quantization step size determined by the bit depth and prescribed range. Prior to quantization, the input is scaled to a normalized range and then the quantized result is rescaled to its original range. Given that the DAC and ADC may have different bit depths, we denote their respective quantization operators by $Q_{\text{DAC}}$ and $Q_{\text{ADC}}$. Extending (\ref{eq:ofterrs}), the complete simulation model for the OFT can be expressed as
\begin{equation} \label{eq:2f_errmodel}
    Q_{\text{ADC}}\left(\mathcal{H}\odot F\left\{Q_{\text{DAC}}(T)\right\}+\varepsilon\right),
\end{equation}
where now both optical errors and data quantization are included in the model.

\subsection{OFT-Based Preconditioner} \label{sec:opt_implementation}

\begin{figure}[h]
     \centering
         \centering
         \includegraphics[width=.5\textwidth]{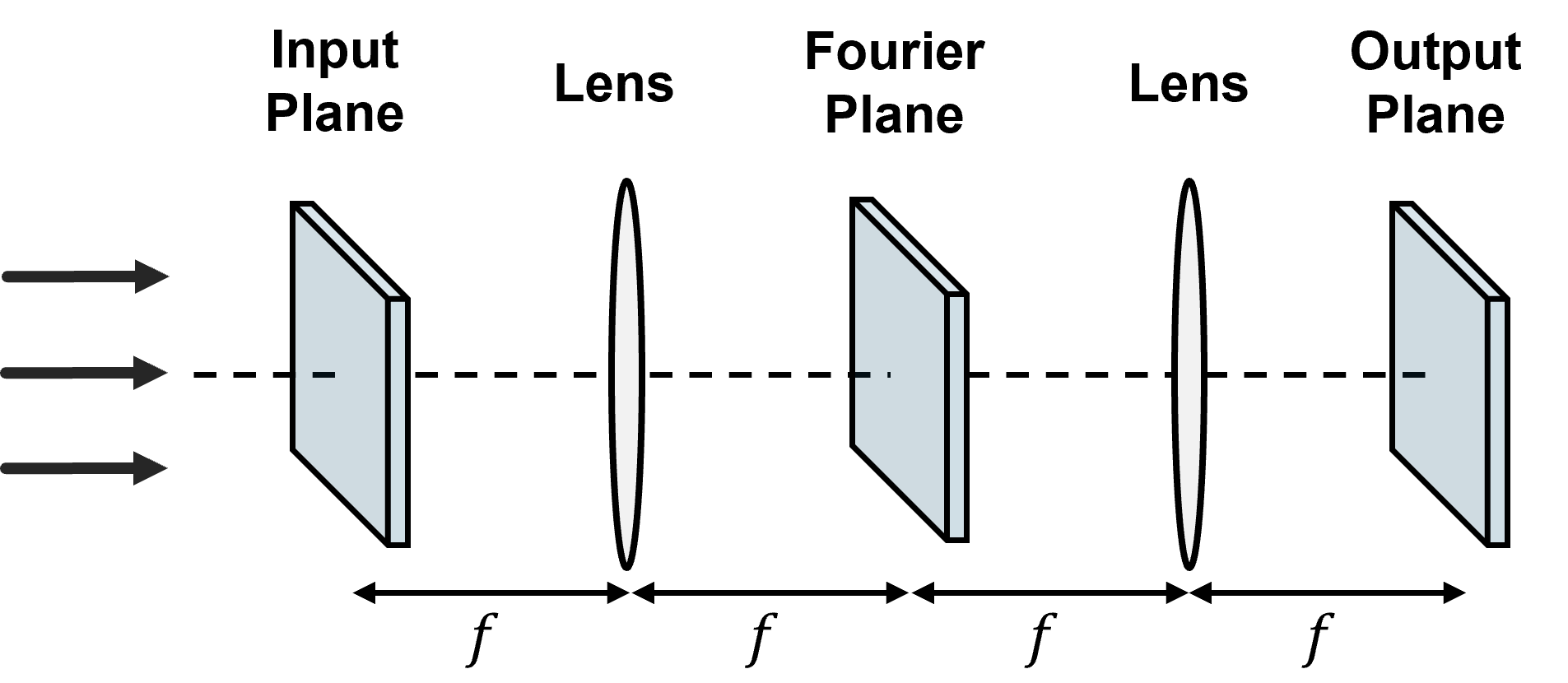}
         \caption{Schematic of a $4f$ optical system where each component is separated by focal length $f$. The first lens performs the Fourier transform and the second lens performs the inverse Fourier transform. A programmable optical light modulator, when placed in the Fourier plane, performs element-wise multiplication in the Fourier domain. }
         \label{fig:4fconfig}
\end{figure}

Recall the two types of preconditioner considered in this work: Type 1 (\ref{eq:type1precon}) which may use a Fourier diagonalizable approximation of $A$ and Type 2 (\ref{eq:type2precon}) which uses the LDL block factorization with a Fourier diagonalizable approximation of $(A^TA+\alpha^2 I)$. For the Type 1 preconditioner and its necessary Fourier transform operations, we consider the single-lens 2$f$ optical system, as discussed in Section~\ref{sec:oft}. For the Type 2 preconditioner, however, its application, $s=M^{-1}y$, is given by
\begin{equation*}
\begin{aligned}
\begin{bmatrix}
s_1\\
s_2
\end{bmatrix}
& =
\begin{bmatrix}
I & -A\\
0 & I
\end{bmatrix}
\begin{bmatrix}
I & 0\\
0 & -P^{-1}
\end{bmatrix}
\begin{bmatrix}
I & 0\\
-A^T & I
\end{bmatrix}
\begin{bmatrix}
y_1\\
y_2
\end{bmatrix}\\[0.5em]
& = \begin{bmatrix}
    y_1-AP^{-1}(A^Ty_1-y_2)\\
    P^{-1}(A^Ty_1-y_2)
    \end{bmatrix}\\[0.5em]
    & = \begin{bmatrix}
    y_1 -As_2\\
    P^{-1}(A^Ty_1-y_2)
    \end{bmatrix}
\end{aligned}
\end{equation*}
where notably the expansion of $s_2$ is
\begin{equation*}
s_2 = {F}^{-1}\left\{  \left(\bar{\Lambda}+\alpha^2I\right)^{-1} \,\, {F}\left\{A^Ty_1-y_2 \right\}\right\},
\end{equation*}
where $\bar{\Lambda}$ contains the eigenvalues of a Fourier diagonalizable approximation to $A^TA$. This structure is well-suited for a $4f$ (or two-lens) optical system where a Fourier transform, element-wise multiplication in the Fourier domain, and an inverse Fourier transform can be performed within a single optical system. This well-established optical configuration utilizes an SLM for the input, a lens to perform the Fourier transform, another SLM for the Fourier plane filter, and a secondary lens to perform the inverse Fourier transform, as depicted in Figure~\ref{fig:4fconfig}. For further illustration, specifically for photographs of $4f$ optical systems built in laboratory settings, we refer the reader to \cite{Chang_OpticalNN_2018,Fitzpatrick}.

This system is advantageous because both Fourier transforms and Fourier-domain multiplication can be implemented within the optical system, rather than returning the data to the digital processor, which reduces the number of necessary data conversions and associated conversion time and quantization errors. The amount of detector noise introduced into the preconditioner is also reduced due to the elimination of the intermediate detection stage that would be necessary if the same operations were implemented with a $2f$ system. Accounting for these differences, we can construct a model for the $4f$ system that incorporates both optical and quantization errors, similar to the $2f$ system model (\ref{eq:2f_errmodel}):
\begin{equation}
    Q_{\text{ADC}}\left(
    {F}^{-1}\{\mathcal{H}_2\odot Q_{\text{DAC}}(\bar{\Lambda}) \odot \left(
    \mathcal{H}_1\odot {F}\left\{Q_{\text{DAC}}(T)\right\}+\varepsilon_1
    \right)\}+\varepsilon_2
    \right).
\end{equation}
Here, $\mathcal{H}_1$ and $\mathcal{H}_2$ represent the transfer functions for the forward and inverse transforms, respectively, while $\varepsilon_1$ and $\varepsilon_2$ represent the corresponding additive noise. 

Recall that $\mathcal{H}$ incorporates the effects of finite aperture, aberrations, and lateral misalignment. The aberrations and misalignment for the forward and inverse transforms in practical optical systems will vary and therefore require their own respective transfer functions. Similarly, noise accumulated during the forward and inverse transforms will vary in both type and magnitude and is therefore modeled independently of each other. For instance, detector noise only contributes to $\varepsilon_2$ since the detection stage is only at the output of the inverse transform for the $4f$ system. 

Quantization effects have been incorporated at the input, Fourier plane, and output stages. It is worth noting that if the same Fourier-domain filter is used for each application of the preconditioner, then the filter only needs to be loaded onto the SLM once at the beginning of the iterative process, rather than before each application of the optical system. Thus, quantization is only applied once to the Fourier-domain filter when simulating, rather than at each application of the OFT. 

There are, however, components of the simulated optical systems, both $2f$ and $4f$, that should vary with each application as they would in the physical system. Specifically, the noise introduced into the model should be independently generated for each application. In contrast, we assume that the parameters associated with the following are held relatively constant throughout the iterative process: finite aperture, misalignment, and aberrations. The bit depth utilized for quantization is also held fixed, however, the resulting quantization errors will vary between applications. The variations in the noise and quantization errors will cause the OFT-based preconditioner to change from one iteration to the next, motivating the use of the flexible variant of GMRES.

\section{Numerical Results} \label{sec:results}

This section is organized as follows: Section~\ref{sec:oftparams} details the setup and parameters used for the OFT simulations in the subsequent numerical experiments. Section~\ref{subsec:oftacc} evaluates the simulated OFT independently of its use within the proposed hybrid digital-optical preconditioning framework. Then, in Sections~\ref{subsec:prblurspeck} and~\ref{subsec:prtomo}, we evaluate the effectiveness of the proposed simulated OFT-based preconditioner within Krylov solvers for Tikhonov-regularized inverse problems. We consider GMRES, and its flexible variant FGMRES for our Krylov solvers and consider the following two test problems in Sections~\ref{subsec:prblurspeck} and~\ref{subsec:prtomo}, respectively: \emph{Speckle} and \emph{Tomo} from the software packages IR tools \cite{IRtools} and AIR tools II \cite{AIRTools}. \emph{Speckle} simulates spatially invariant blurring with a speckle point spread function (PSF) that models atmospheric turbulence and \emph{Tomo} generates X-ray tomographic reconstruction problems based on the Radon transform.

The reported numerical results of this work were carried out in MATLAB R2022b 64-bit on a Lenovo laptop running Microsoft Windows 11 Pro with a 12th Gen Intel Core i7-1255U processor and 16 GB of RAM. Results reported herein utilize double precision with unit round-off error of approximately $2.22\times10^{-16}$.

\subsection{OFT Simulation Parameters} \label{sec:oftparams}
For all numerical experiments presented in this work, we consider three optical system quality levels for the simulated OFT: near-ideal, moderate, and degraded. Each quality level defines a range of values for the optical system parameters, as shown in Table~\ref{tab:optparams}. We implement parameter ranges to account for the potential variations in optical system components and conditions in practical implementations. This allows us to evaluate the robustness of both the simulated OFT and the resulting preconditioner across varying levels of optical system quality. For each error source incorporated into the simulated OFT, we define the parameters to be used in our numerical experiments. 

For the effects of finite aperture, we consider the circular pupil function described by (\ref{eq:pupil}) and specify a percentage that relates the radius $r$ to the largest circle inscribed within the square discrete Fourier grid. A value of $100\%$ corresponds to the largest possible pupil radius, while lower percentages increasingly truncate higher spatial frequencies. For lateral misalignment, the displacement of the input relative to the optical axis, in both x- and y-directions, is incorporated according to (\ref{eq:shift}). The magnitude of the misalignment is controlled by the input shift parameter and specified in pixels.

As discussed in Section~\ref{sec:lens_abberations}, optical aberrations can be incorporated through Zernike polynomials (\ref{eq:zernikepoly}). For our simulation purposes, we include the following aberrations: defocus, $x$- and $y$-primary astigmatism, $x$- and $y$-coma, $x$- and $y$-trefoil, and primary spherical; see \cite{Schmidt,Niu_Zernike} for further discussion of these components of the Zernike polynomial. The coefficients of the Zernike terms determine the severity of the corresponding aberration, as well as the overall wavefront error for the system. For our simulation purposes, the Zernike term coefficients are determined based on a prescribed system root-mean-square wavefront error (RMS WFE), which serves as a quantitative measure of wavefront aberration. Specifically, RMS WFE quantifies the overall magnitude of wavefront aberration across the pupil. For orthonormal Zernike polynomials with the piston (uniform phase) term excluded, the coefficients are directly related to the system RMS WFE by
\begin{equation}
    \left\langle W^2(\rho,\theta) \right\rangle^{1/2} = \left(\sum_j a_j^2\right)^{1/2},
\end{equation}
where $W(\rho,\theta)$ is the wavefront error, $\langle \cdot \rangle$ denotes averaging over the pupil, and $a_j$ denotes the coefficient of the $j$th Zernike term \cite{Mahajan_Zernike}. We assign relative weights to the Zernike coefficients and then scale according to the RMS WFE randomly selected from the ranges in Table \ref{tab:optparams}. The near-ideal system RMS WFE range was chosen to represent diffraction-limited performance, while the moderate and degraded ranges correspond to increased wavefront aberration.

The magnitude of the additive complex Gaussian noise is prescribed by a signal-to-noise ratio (SNR) in decibels (dB),
\begin{equation*}
\mathrm{SNR}
=
20\log_{10}
\left(
\frac{\|S\|_F}
{\|\varepsilon\|_F}
\right),
\end{equation*}
where $\|\cdot\|_F$ denotes the Frobenius norm, $\varepsilon$ is the complex Gaussian noise, and $S$ is the forward or inverse transform of the input. Three different SNR parameters are included in Table \ref{tab:optparams}. The `$2f$ forward OFT SNR' and the `$4f$ forward OFT SNR' represent the SNRs used for the forward transform processes for the $2f$ and $4f$ systems, respectively. A larger SNR is chosen for the $4f$ forward transform process due to the absence of detector noise at this stage of the $4f$ system. The `$2f/4f$ inverse OFT SNR' parameter represents the SNR used for the inverse transform process for both optical systems. Lastly, for quantization, the DAC and ADC bit depth parameters determine the resolution of the digital-to-analog and analog-to-digital conversion, respectively. As shown in Table \ref{tab:optparams}, we use fixed bit depths across all optical system quality levels. 

\begin{table}[h] 
\centering \small
\renewcommand{\arraystretch}{1.3}
\setlength{\tabcolsep}{4pt}
\caption{Simulated OFT parameter value ranges considered for optical system quality levels of near-ideal, moderate, and degraded.  }
\begin{tabular}{c | l | c | c | c}
\textbf{Error Source} &
\textbf{Parameter} &
\textbf{Near-Ideal} &
\textbf{Moderate} &
\textbf{Degraded} \\
\Xhline{1.5pt} 
finite aperture
& pupil radius (\%) 
& $90$ - $100$
& $80$ - $90$
& $70$ - $80$ \\ 
misalignment
& input shift (pixels)
& $0.0$ - $0.1$
& $0.1$ - $0.5$
& $0.5$ - $1.0$ \\ 
aberrations
& system RMS WFE (waves)
& $0.00$ - $0.05$
& $0.05$ - $0.10$
& $0.10$ - $0.15$ \\ \hline
\multirow{3}{*}{noise}
& $2f$ forward OFT SNR (dB)
& $45$ - $50$
& $35$ - $45$
& $25$ - $35$ \\
& $4f$ forward OFT SNR (dB)
& $55$ - $60$
& $45$ - $55$
& $35$ - $45$ \\
& $2f$/$4f$ inverse OFT SNR (dB)
& $45$ - $50$
& $35$ - $45$
& $25$ - $35$ \\
\hline
\multirow{2}{*}{quantization}
&DAC bit depth & \multicolumn{3}{c}{8 bits} \\
&ADC bit depth & \multicolumn{3}{c}{12 bits} \\
\hline
\end{tabular}
\label{tab:optparams}
\end{table}

\subsection{Simulated OFT Accuracy} \label{subsec:oftacc}

We first evaluate the simulated OFT independently of its use within a preconditioner. In this section we examine the effects of the optical error sources and quantization on the Fourier operator through its singular values and then we evaluate the accuracy of the simulated OFT when applied directly to image data. 

It is known that the normalized Fourier transform is a unitary operator whose singular values are equal to $1$. Thus, for the optical Fourier transform, we investigate whether its modeled errors cause a deviation from the ideal unitary behavior, and if so, how much deviation is present and which singular values are affected. To investigate this, we compare the singular values for four different normalized Fourier operators: the digital discrete Fourier transform and the simulated OFT for all three optical quality levels presented in Table \ref{tab:optparams}.
Considering a $64 \times 64$ spatial grid, we construct the matrix representation of each operator by applying the respective transform to $4096$ standard basis vectors, where the outputs form the columns of the resulting $4096\times4096$ matrix. The digital Fourier transform is computed using MATLAB's \texttt{fft2} and denoted by FFT. The simulated optical Fourier transforms are computed according to the $2f$ system error model (\ref{eq:2f_errmodel}). For each optical quality, the parameters for each error source are randomly selected from the range of values provided in Table \ref{tab:optparams}.

\begin{figure}[h]
\centering
\includegraphics[width=.5\linewidth]{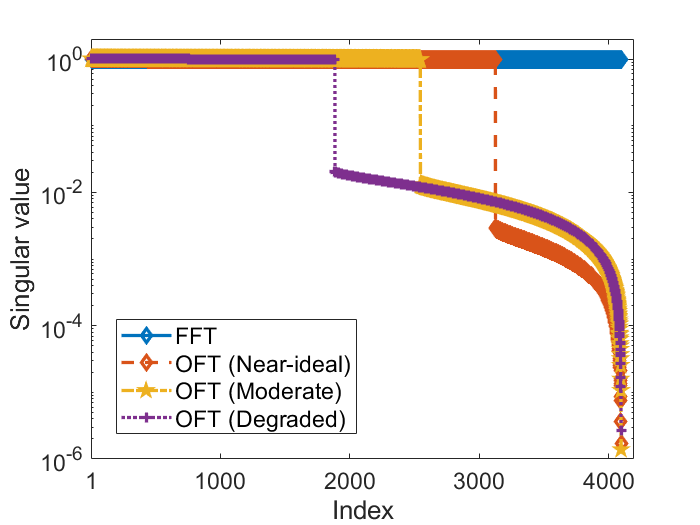}
\caption{Comparison of singular values of the normalized digital Fourier transform (FFT) operator and normalized simulated optical Fourier transform (OFT) operators with near-ideal, moderate, and degraded optical system qualities. A $64\times64$ spatial grid was considered, corresponding to Fourier operators of size $4096\times4096$.}
\label{fig:singvals_comparison}
\end{figure}

In Figure~\ref{fig:singvals_comparison}, we present the singular values of each operator in descending order. The normalized digital FFT has $4096$ singular values numerically equal to $1$ in double precision (i.e., are within double precision round-off error), as expected. The simulated OFT operators, however, deviate from this unitary behavior. Specifically, the singular values align with the digital FFT singular values for the first approximately $1900-3100$ values, but then drop sharply by roughly $2$ orders of magnitude and decay towards $10^{-6}$ for the remaining singular values. The index at which the deviation from $1$ occurs is dependent on the optical quality. With decreasing quality (increased optical system errors), the number of singular values aligned with the digital FFT decreases, with approximately $3100$ aligned for the near-ideal case and approximately $1900$ aligned for the degraded case. Nevertheless, despite the errors introduced by the simulated optical system, a significant portion of the singular values of the digital Fourier operator are preserved. 

We now directly compare the digital and optical Fourier transforms of image data. For this comparison, we use the Hubble space telescope image from the IR tools \cite{IRtools} software package and incorporate varying levels of additive Gaussian noise into the image data. Specifically, the noisy image is given by $x=x_{true}+e$, where $e$ consists of mean zero normally distributed random entries that are scaled to correspond to a prescribed noise level given by $100(\|e\|/\|x_{\text{true}}\|)\%$. We consider noise levels of $1\%, \, 5\%, \text{and } 10\%$.

In Table \ref{tab:PRblur_2f4ferror}, we report the following relative errors: the error between the simulated OFT ($2f$ system) applied to image $x$ and the digital FFT applied to image $x$, the error after applying both the forward and inverse simulated OFT to the image $x$ with both the $2f$ and $4f$ systems, and the error after applying both the forward and inverse digital FFT to the image $x$. Errors involving the simulated OFT were computed across $100$ independent trials and the mean and standard deviation of the relative errors are reported. For each trial, the optical system parameters are randomly selected from the appropriate optical system quality range, given by Table \ref{tab:optparams}. 

MATLAB's \texttt{fft2} and \texttt{ifft2} are able to recover $x$ within the round-off error of double precision, as shown in the final column of the table. The simulated optical transforms, however, introduce nonneglible errors into the reconstruction. As expected, we observe decreased accuracy with lower optical system quality, or equivalently, with increased optical errors. We similarly observe decreased accuracy, although less pronounced, with increased noise in $x$. The relative errors produced by the $4f$ system are noticeably smaller than those produced by the $2f$ system. This is due to the reduction in detection-related errors (detector noise and ADC quantization) in the $4f$ implementation. Thus, the $4f$ configuration should be utilized when feasible. 

The results shown in Table \ref{tab:optparams} demonstrate that the errors introduced by the simulated OFTs are nonnegligible across all optical system quality levels and system configurations. These results not only provide a comparison between the digital and simulated optical transforms, but they also provide some insight into the amount of error that will be introduced into the preconditioner when using OFTs. 

\begin{table}[h]
\footnotesize\centering
\renewcommand{\arraystretch}{1.5}
\setlength{\tabcolsep}{4pt}
\caption{Accuracy of the simulated OFT, relative to the digital FFT, using the Hubble space telescope image, with $1\%$, $5\%$, and $10\%$ noise levels, as input $x$. The OFT and digital FFT are denoted by ${F}_{\mathrm{opt}}$ and ${F}$, respectively. The relative errors for the forward OFT and the combined forward-inverse OFT are shown for the three optical quality levels considered: near-ideal, moderate, and degraded. For each quality level, the OFT-based results are computed over $100$ independent trials with the mean and standard deviation reported. The relative error after both forward and inverse digital FFT is included for reference.}
\begin{tabular}{c | c | c | c | c | c}
\multirow{2}{*}{\makecell{{Noise}\\{in ${x}$}}} &
\multirow{2}{*}{\makecell{{Optical}\\{Quality}}} &
\multirow{2}{*}{$\displaystyle
\frac{\|{F}_{\mathrm{opt}}(x)-{F}(x)\|}
     {\|{F}(x)\|}$} &
\multicolumn{2}{c|}{$\displaystyle
\frac{\|{F}_{\mathrm{opt}}^{-1}
({F}_{\mathrm{opt}}(x))-x\|}
{\|x\|}$} &
\multirow{2}{*}{$\displaystyle
\frac{\|{F}^{-1}({F}(x))-x\|}
     {\|x\|}$} \\[0.5em]
& & & $2f$ system & $4f$ system& \\

\Xhline{1.5pt} 
\multirow{3}{*}{$1\%$} & Near-ideal & $0.1640 \pm 0.1171$ & $0.1672 \pm 0.0711$ & $0.0702 \pm 0.0196$ &\multirow{3}{*}{$3.54e-16$} \\
& Moderate & $0.4370 \pm 0.2545$ & $0.5000 \pm 0.3609$ & $0.1804 \pm 0.1056$ &\\
& Degraded & $0.7064 \pm 0.3828$ & $0.9413 \pm 0.6791$ & $0.3674 \pm 0.2419$ &\\
\hline \hline

\multirow{3}{*}{$5\%$} & Near-ideal & $0.1672 \pm 0.1155$ & $0.1718 \pm 0.0700$ & $0.0755 \pm 0.0191$ &\multirow{3}{*}{$3.65e-16$} \\
& Moderate & $0.4393 \pm 0.2527$ & $0.5029 \pm 0.3584$ & $0.1853 \pm 0.1042$ &\\
& Degraded & $0.7084 \pm 0.3809$ & $0.9433 \pm 0.6760$ & $0.3714 \pm 0.2401$ &\\
\hline \hline

\multirow{3}{*}{$10\%$} & Near-ideal & $0.1761 \pm 0.1112$ & $0.1837 \pm 0.0670$ & $0.0895 \pm 0.0180$ &\multirow{3}{*}{$3.91e-16$} \\
& Moderate & $0.4461 \pm 0.2476$ & $0.5113 \pm 0.3514$ & $0.1990 \pm 0.1004$ &\\
& Degraded & $0.7141 \pm 0.3752$ & $0.9489 \pm 0.6673$ & $0.3831 \pm 0.2350$ &\\
\hline
\end{tabular}
\label{tab:PRblur_2f4ferror}
\end{table}

\subsection{Image Deblurring} \label{subsec:prblurspeck}

\begin{figure}[h]
\centering
\begin{subfigure}{0.3\textwidth}
    \centering
    \includegraphics[width=\linewidth,height=3.5cm,keepaspectratio]{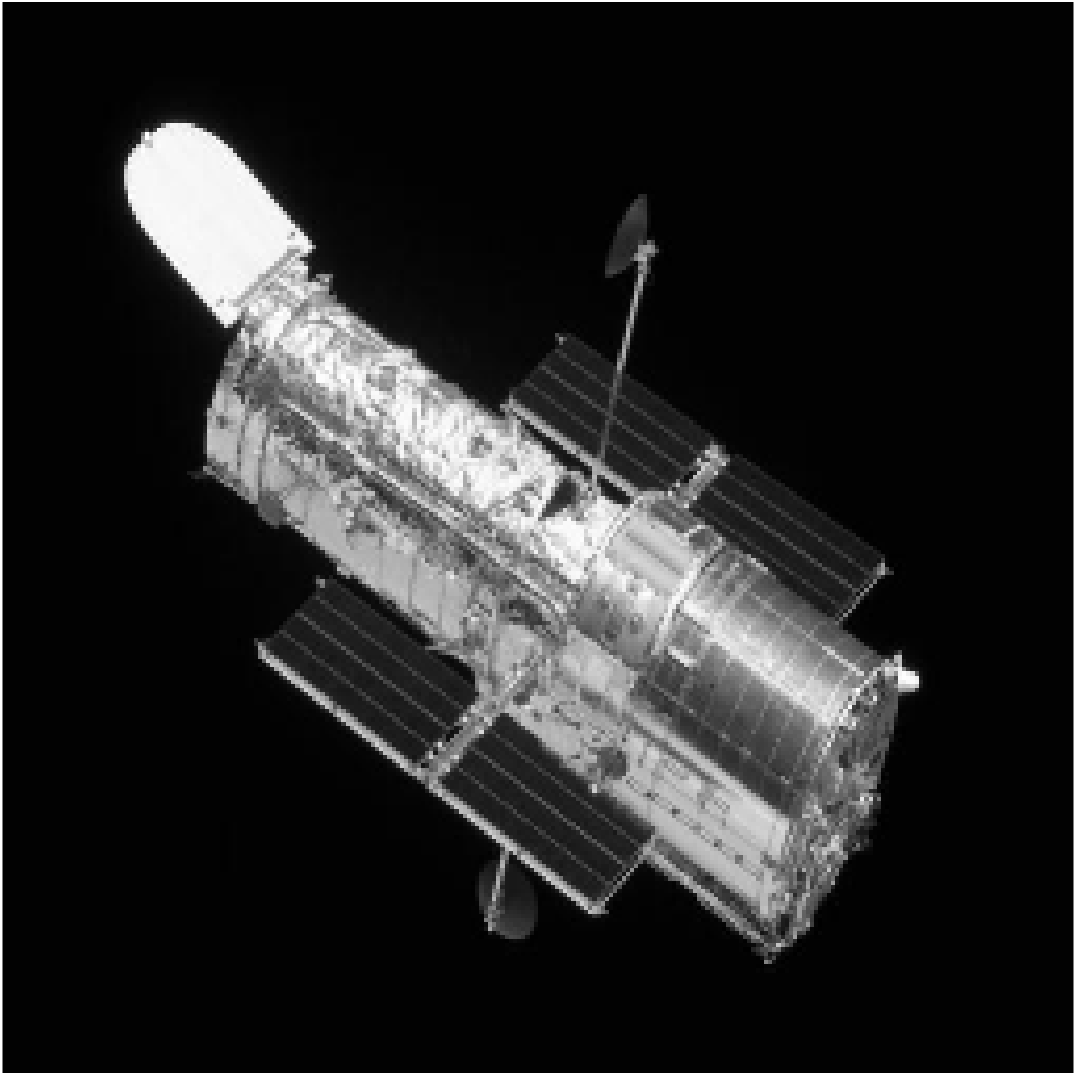}
    \caption{}
    \label{fig:prblurspeck_xtrue}
\end{subfigure}
\begin{subfigure}{0.3\textwidth}
    \centering
    \includegraphics[width=\linewidth,height=3.5cm,keepaspectratio]{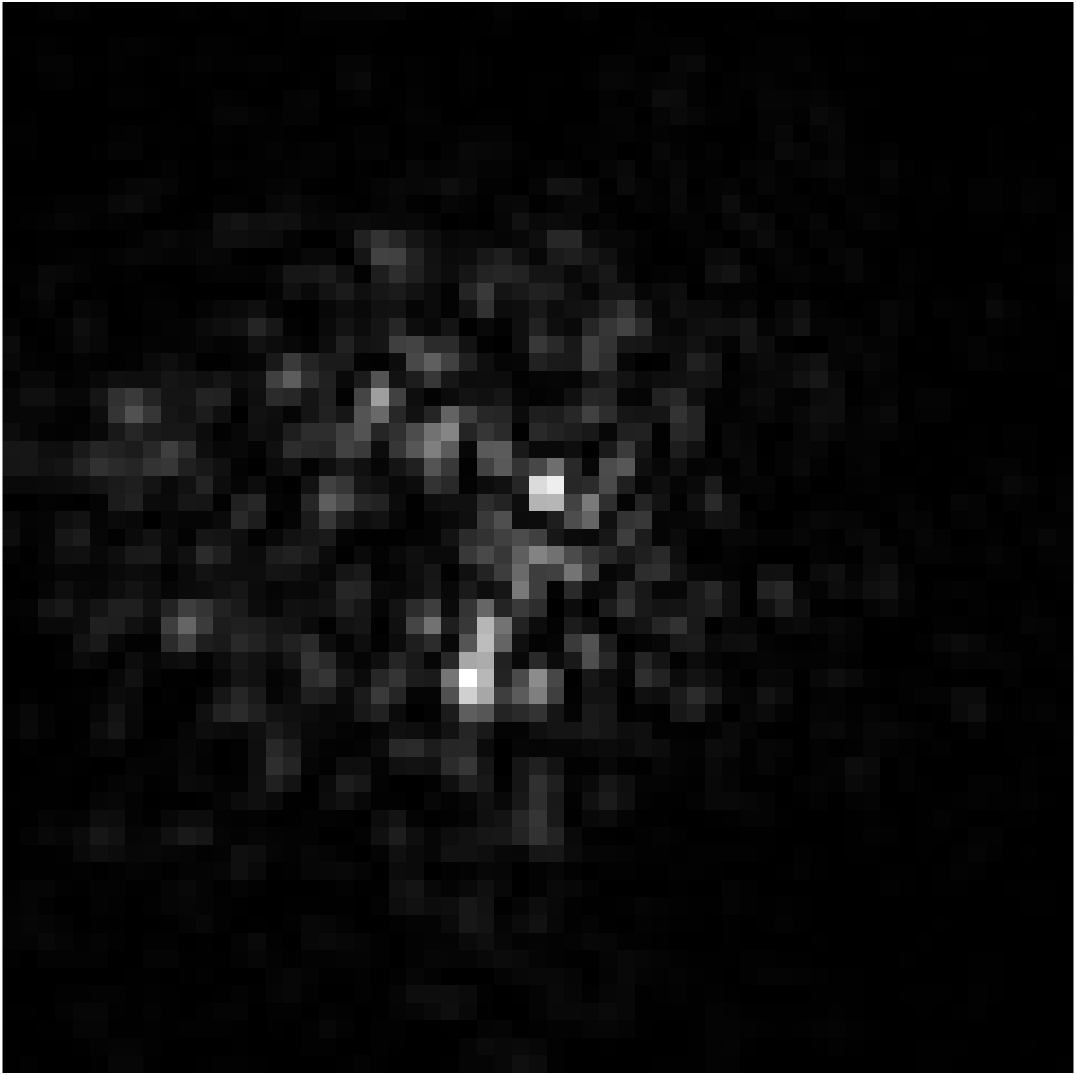}
    \caption{}
    \label{fig:prblurspeck_psf}
\end{subfigure}
\begin{subfigure}{0.3\textwidth}
    \centering
    \includegraphics[width=\linewidth,height=3.5cm,keepaspectratio]{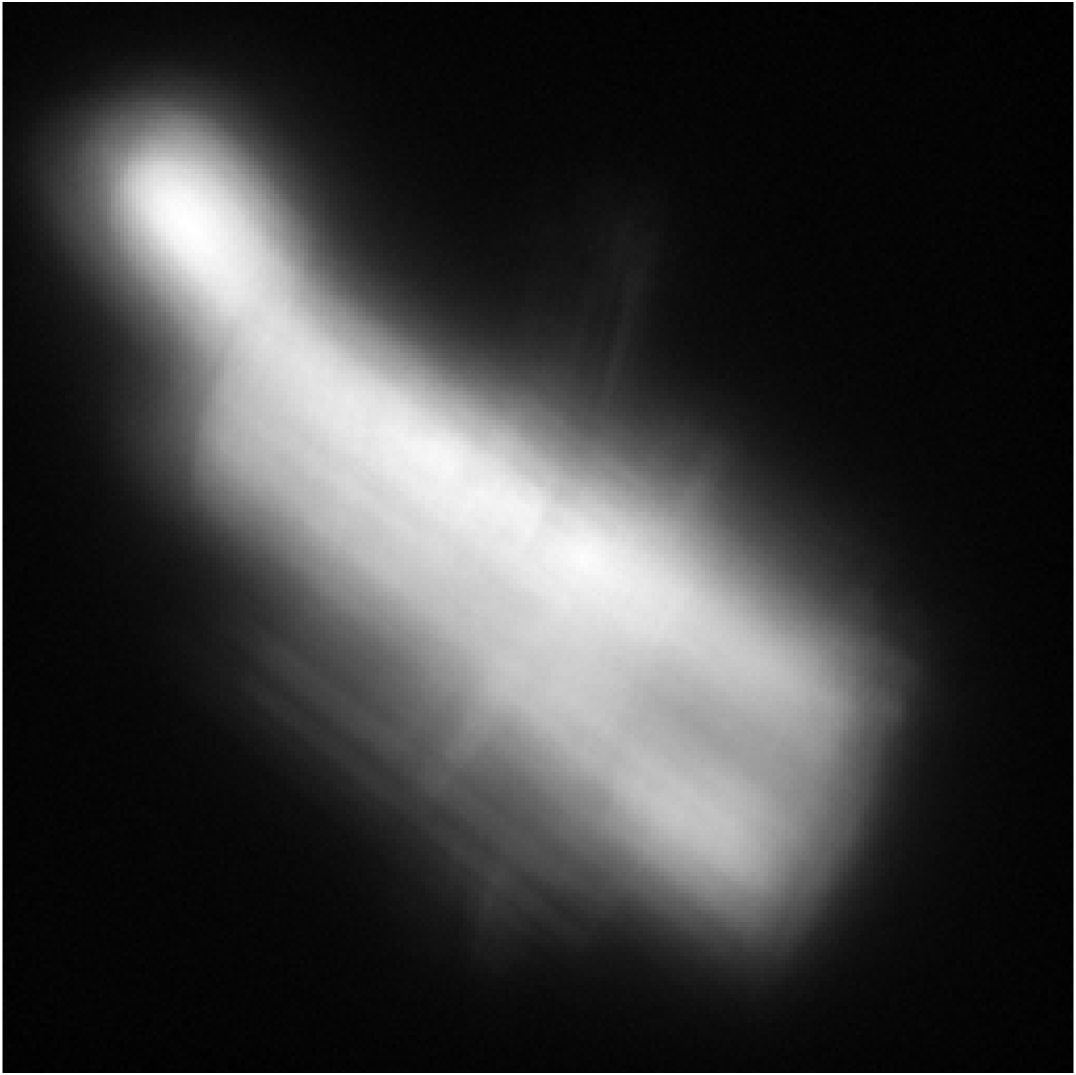}
    \caption{}
    \label{fig:prblurspeck_bn}
\end{subfigure}
\vspace{-0.8em}
\caption{$2$D image deblurring problem \emph{Speckle}: (a) true Hubble image ($256\times 256$ pixels), (b) zoomed-in view of the \textit{speckle} point spread function (central $61\times61$ pixels), and (c) blurred Hubble image with $1\%$ additive Gaussian noise.}
\label{fig:prblurspecksetup}
\end{figure}

We now demonstrate the effectiveness of the proposed simulated OFT-based preconditioner within GMRES and FGMRES for the \emph{Speckle} $2$D deblurring test problem. In this test problem, the forward operator $A$ represents convolution of the true image with a \emph{Speckle} PSF that models spatially invariant atmospheric turbulence with imposed zero boundary conditions. We consider images with $256\times256$ pixels, corresponding to matrix $A$ of size $65536\times65536$. Additive Gaussian noise is introduced to the right-hand side $b=Ax+e$, where $e$ consists of mean zero normally distributed random entries that are scaled to correspond to a prescribed noise level given by $100(\|e\|/\|Ax_{\text{true}}\|)\%$, where $x_{\text{true}}$ denotes the true vectorized image. We consider noise levels of $1\%, \, 5\%, \text{and } 10\%$. The true image and \emph{Speckle} PSF are shown in Figures~\ref{fig:prblurspeck_xtrue} and ~\ref{fig:prblurspeck_psf}, respectively. The corresponding blurred image with $1\%$ additive Gaussian noise is shown in Figure~\ref{fig:prblurspeck_bn}.

To reconstruct the image from the noisy, blurred data, we approximate the solution to the saddle-point formulation shown in (\ref{eq:saddlePoint}) using GMRES and FGMRES with right preconditioning, as described in (\ref{eq:min_saddleProb}). To model long-exposure atmospheric blur, a Gaussian PSF is frequently adopted as a smooth isotropic approximation to the atmospheric turbulence model. If we approximate the zero boundary condition of our underlying operator with a periodic boundary condition, which is common when the image of interest stems from astronomical imaging, then we may approximate the blurring operator by a block circulant with circulant blocks matrix which may be diagonalized by the FFT. Thus, for this problem, we consider the Type 1 preconditioner introduced in (\ref{eq:type1precon}) that utilizes a single-lens 2f optical system setup. 

For the regularization parameter, we select an ``optimal'' $\alpha$ with respect to the error norm and hold this value fixed throughout the iterations. We use the following minimizer
\begin{equation} \label{eq:regparam}
    \alpha = \arg \min_\alpha \|\widetilde{x}_{\alpha} - x_{\text{true}}\|,
\end{equation}
where $\widetilde{x}_{\alpha}$ denotes the solution of the nearby approximate Tikhonov problem
\begin{equation*}
    \widetilde{x}_{\alpha}=\arg\min_x\left\{\|\bar{C}x-b\|^2+\alpha^2\|x\|^2\right\}
\end{equation*}
obtained by replacing $A$ with its circulant approximation $\bar{C}$.
We acknowledge this selection method is of course not a reasonable one in practice due to its dependence on knowledge of the true solution. However, our focus is on evaluating the performance of the proposed OFT-based preconditioners regardless of the regularization parameter selection. We note that there exist many regularization parameter choice methods for Tikhonov regularization problems such as the discrepancy principle, generalized cross validation, or the L-curve criterion; see e.g., \cite{ChungGazzola_survey, HNO, RR13} for further discussion.

The discrepancy principle is utilized to determine a suitable stopping criterion. This technique stops an iterative solver once the residual norm of the current solution drops to or below a prescribed threshold based on available knowldege of the noise contaminating the data. Specifically, the discrepancy principle terminates the iterative method at the first iteration $p$ for which
\begin{equation} \label{eq:discrepancy}
    \|Ax_{p}-b \|\leq\tau \delta
\end{equation}
where the constant $\tau >1$ is considered as a ``safety factor'' and $\delta$ is an upper bound on the norm of the noise contaminating the data, such that $\|b-b_{\text{true}}\|=\|e\|\leq \delta$. We set $\tau = 1.05$ for all numerical experiments. Although GMRES and FGMRES are applied to the saddle-point formulation, the discrepancy principle is evaluated using the residual of the available problem \eqref{leastSquares}.

To evaluate the quality of the reconstructed images, we compute the relative reconstruction error (RRE) defined by
\begin{equation} \label{eq:RRE}
\text{RRE}\left(x_{p}\right)=\frac{\left\| x_{p}-x_{\text{true}}\right\|}{\left\|x_{\text{true}}\right\|},
\end{equation}
where $x_{p}$ is the reconstructed image at the $p^{\text{th}}$ iterate satisfying the stopping criterion.

In the simulated OFTs, we consider the $2f$ optical system with all three optical quality levels: near-ideal, moderate, and degraded. For each method utilizing an OFT-based preconditioner, $50$ independent trials are run, where for each trial the parameters for the optical error sources are randomly sampled from the corresponding range provided in Table \ref{tab:optparams}. The parameters associated with finite aperture, misalignment, and aberrations are held fixed across iterations for a given trial. The parameters associated with noise in the optical system are varied across iterations, meaning at each iteration the noise is randomly generated using an SNR value randomly sampled from the prescribed range. The quantization parameters are held fixed across optical qualities and iterations, but the resulting quantization errors vary with each application of the quantization operator.

We consider the following methods for comparison: GMRES with no preconditioning, digital FFT-based preconditioning, and simulated OFT-based preconditioning, as well as FGMRES with simulated OFT-based preconditioning. In Table \ref{tab:prblurspeck_results_M1}, we report the RRE for each of these methods across varying noise levels in $b$.  The reported RRE values provided correspond to the $p^{\text{th}}$ iterative solution when the discrepancy principle is satisfied, or alternatively, when a maximum number of iterations of $100$ is achieved. For the OFT-based methods, we consider three different levels of optical system quality, as described by Table \ref{tab:optparams}. For each level, the mean and standard deviation ($\mu \pm \sigma$) of the $50$ RREs and corresponding breakout iterations are reported. We also report the minimum and maximum RRE and their corresponding breakout iterations.

In Table \ref{tab:prblurspeck_results_M1}, we observe increased reconstruction error as optical system quality decreases across all OFT-based methods. This effect of optical system quality appears more pronounced for the lower noise levels in $b$, particularly for the $1\%$ noise case. Most notably, we observe that FGMRES with OFT-based preconditioning, has comparable performance to the digital FFT-based version, in terms of both reconstruction error and iteration count for the near-ideal and moderate optical quality ranges across all noise levels considered in $b$. In contrast, the non-flexible version of GMRES with OFT-based preconditioning exhibits larger reconstruction errors for the cases of $1\%$ and $5\%$ noise in $b$, in comparison to the digital FFT-based method. These results suggest that given sufficient optical system quality and the use of a flexible Krylov solver, OFT-based preconditioners can achieve reconstruction quality and convergence speedup comparable to that of digitally computed FFT-based preconditioners.

\begin{table}[h]
\centering \footnotesize
\renewcommand{\arraystretch}{1.3}
\setlength{\tabcolsep}{4pt}
\caption{RRE results for the \emph{Spectra} test problem with $1\%$, $5\%$, and $10\%$ noise. GMRES with no preconditioning, GMRES with digital FFT-based preconditioning, GMRES with OFT-based preconditioning, and FGMRES with OFT-based preconditioning are considered. For the OFT-based methods, the optical system quality is varied among near-ideal, moderate, and degraded. Results for the OFT-based methods are computed over $50$ independent trials. The mean and standard deviation of the RRE and the corresponding iteration count is reported. Also reported for the OFT-based methods is the minimum and maximum RRE with their corresponding iteration counts. The RRE values presented correspond to the $p^{\text{th}}$ iterative solution when the discrepancy principle is satisfied, or alternatively, when the specified maximum number of iterations is achieved ($p=100$).} 
\begin{tabular}{c | l | c | l | l | l}
\multirow{2}{*}{\makecell{\textbf{Noise}\\\textbf{in $\boldsymbol{b}$}}} & \multirow{2}{*}{\textbf{Method}} & \multirow{2}{*}{\makecell{\textbf{Optical}\\\textbf{Quality}}} & \multicolumn{3}{c}{\textbf{RRE}} \\
\cline{4-6}
& & & \textbf{$\boldsymbol{\mu \pm \sigma}$ (Iter. $\boldsymbol{\mu \pm \sigma}$)} &
\textbf{Min (Iter.)} &
\textbf{Max (Iter.)} \\
\Xhline{1.5pt}

\multirow{8}{*}{$1\%$} & GMRES & --- & 0.2081 (73) & \hspace{.6cm}--- & \hspace{.6cm}--- \\
\cline{2-6}
& GMRES (dig.) & --- & 0.2019 (4) & \hspace{.6cm}--- & \hspace{.6cm}--- \\
\cline{2-6}
& \multirow{3}{*}{GMRES (opt.)} & Near-ideal & $0.2726 \pm 0.0063$ ($100 \pm 0$) & 0.2608 (100) & 0.2924 (100) \\
& & Moderate & $0.3074 \pm 0.0425$ ($100 \pm 0$) & 0.2637 (100) & 0.4270 (100) \\
& & Degraded & $0.7747 \pm 0.4267$ ($100 \pm 0$) & 0.3385 (100) & 1.9346 (100) \\
\cline{2-6}
& \multirow{3}{*}{FGMRES (opt.)} & Near-ideal & $0.1998 \pm 0.0004$ ($4 \pm 0$) & 0.1992 (4) & 0.2008 (4) \\
& & Moderate & $0.2117 \pm 0.0145$ ($5 \pm 14$) & 0.1984 (4) & 0.2811 (100) \\
& & Degraded & $0.5645 \pm 0.3432$ ($90 \pm 29$) & 0.2173 (4) & 1.8382 (100) \\
\hline \hline 

\multirow{8}{*}{$5\%$} & GMRES & --- & 0.2995 (29) & \hspace{.6cm}--- & \hspace{.6cm}--- \\
\cline{2-6}
& GMRES (dig.) & --- & 0.2900 (1) & \hspace{.6cm}--- & \hspace{.6cm}--- \\
\cline{2-6}
& \multirow{3}{*}{GMRES (opt.)} & Near-ideal & $0.3296 \pm 0.0013$ ($2 \pm 0$) & 0.3256 (2) & 0.3317 (2) \\
& & Moderate & $0.3144 \pm 0.0058$ ($2 \pm 0$) & 0.3044 (2) & 0.3252 (2) \\
& & Degraded & $0.4132 \pm 0.1391$ ($49 \pm 48$) & 0.3041 (31) & 0.7843 (100) \\
\cline{2-6}
& \multirow{3}{*}{FGMRES (opt.)} & Near-ideal & $0.2897 \pm 0.0005$ ($2 \pm 0$) & 0.2886 (2) & 0.2907 (2) \\
& & Moderate & $0.2926 \pm 0.0051$ ($2 \pm 0$) & 0.2825 (2) & 0.3078 (2) \\
& & Degraded & $0.4024 \pm 0.1392$ ($46 \pm 49$) & 0.2945 (2) & 0.7657 (100) \\
\hline \hline 

\multirow{8}{*}{$10\%$} & GMRES & --- & 0.3415 (17) & \hspace{.6cm}--- & \hspace{.6cm}--- \\
\cline{2-6}
& GMRES (dig.) & --- & 0.3224 (1) & \hspace{.6cm}--- & \hspace{.6cm}--- \\
\cline{2-6}
& \multirow{3}{*}{GMRES (opt.)} & Near-ideal & $0.3250 \pm 0.0007$ ($1 \pm 0$) & 0.3237 (1) & 0.3268 (1) \\
& & Moderate & $0.3294 \pm 0.0067$ ($1 \pm 0$) & 0.3226 (1) & 0.3458 (1) \\
& & Degraded & $0.3785 \pm 0.0777$ ($4 \pm 14$) & 0.3229 (1) & 0.6376 (1) \\
\cline{2-6}
& \multirow{3}{*}{FGMRES (opt.)} & Near-ideal & $0.3250 \pm 0.0007$ ($1 \pm 0$) & 0.3236 (1) & 0.3269 (1) \\
& & Moderate & $0.3295 \pm 0.0067$ ($1 \pm 0$) & 0.3228 (1) & 0.3464 (1) \\
& & Degraded & $0.3799 \pm 0.0803$ ($5 \pm 20$) & 0.3211 (2) & 0.6306 (100) \\
\hline 
\end{tabular}
\label{tab:prblurspeck_results_M1}
\end{table}

\begin{figure}[h]
\centering
\begin{subfigure}[h]{0.32\columnwidth}
    \centering
    \includegraphics[height=3.2cm]{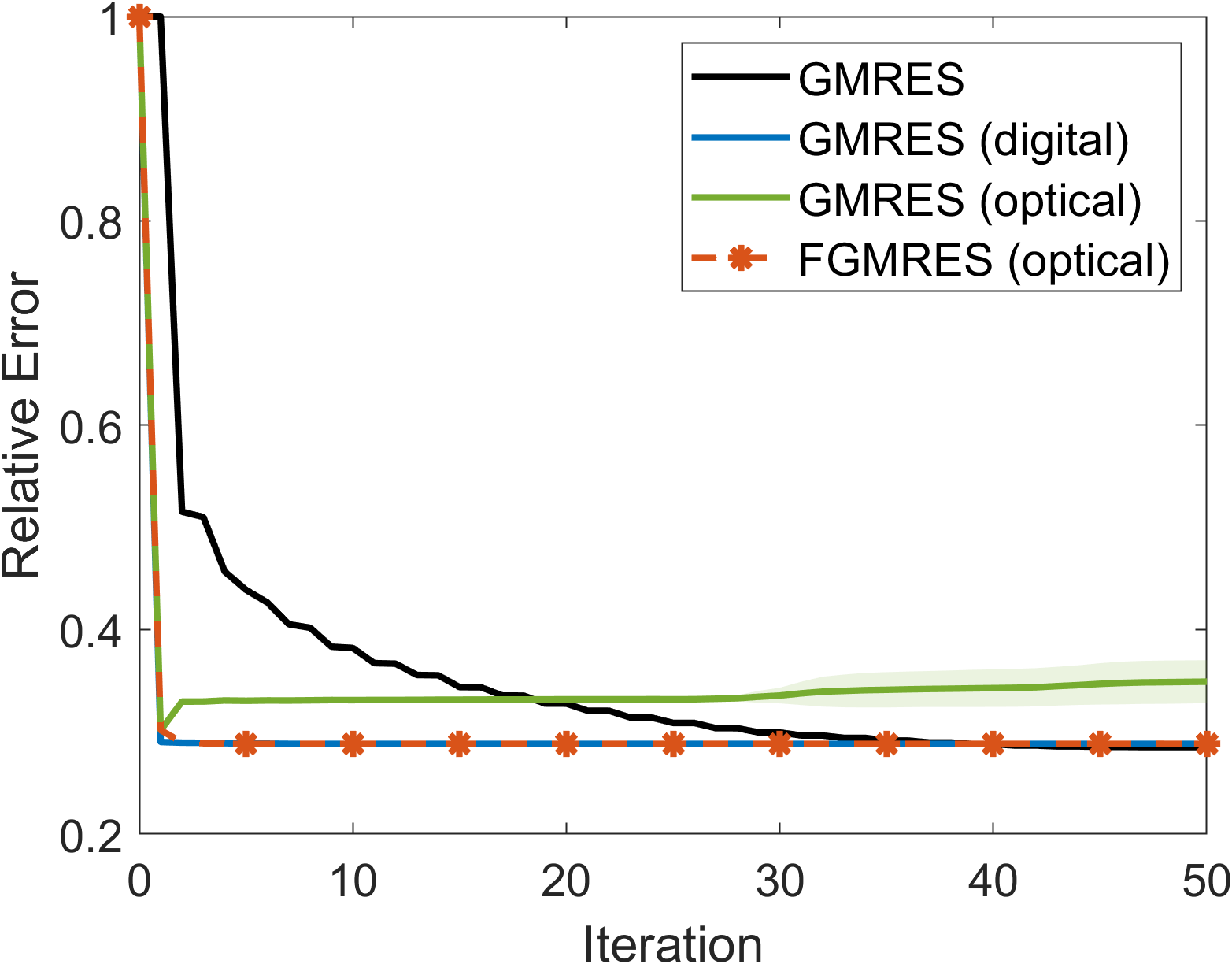}
\end{subfigure}
\hfill
\begin{subfigure}[h]{0.32\columnwidth}
    \centering
    \includegraphics[height=3.2cm]{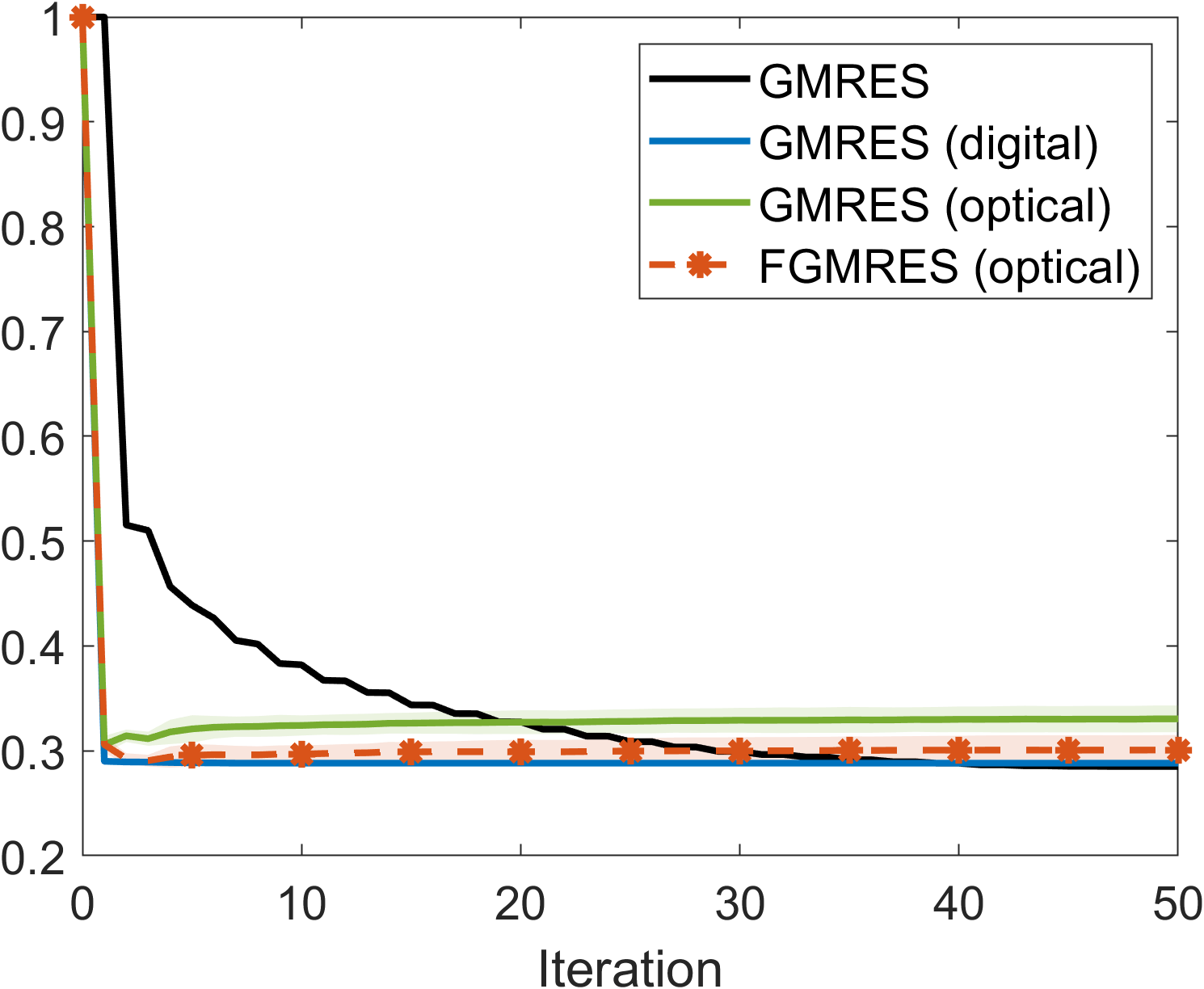}
\end{subfigure}
\hfill
\begin{subfigure}[h]{0.32\columnwidth}
    \centering
    \includegraphics[height=3.2cm]{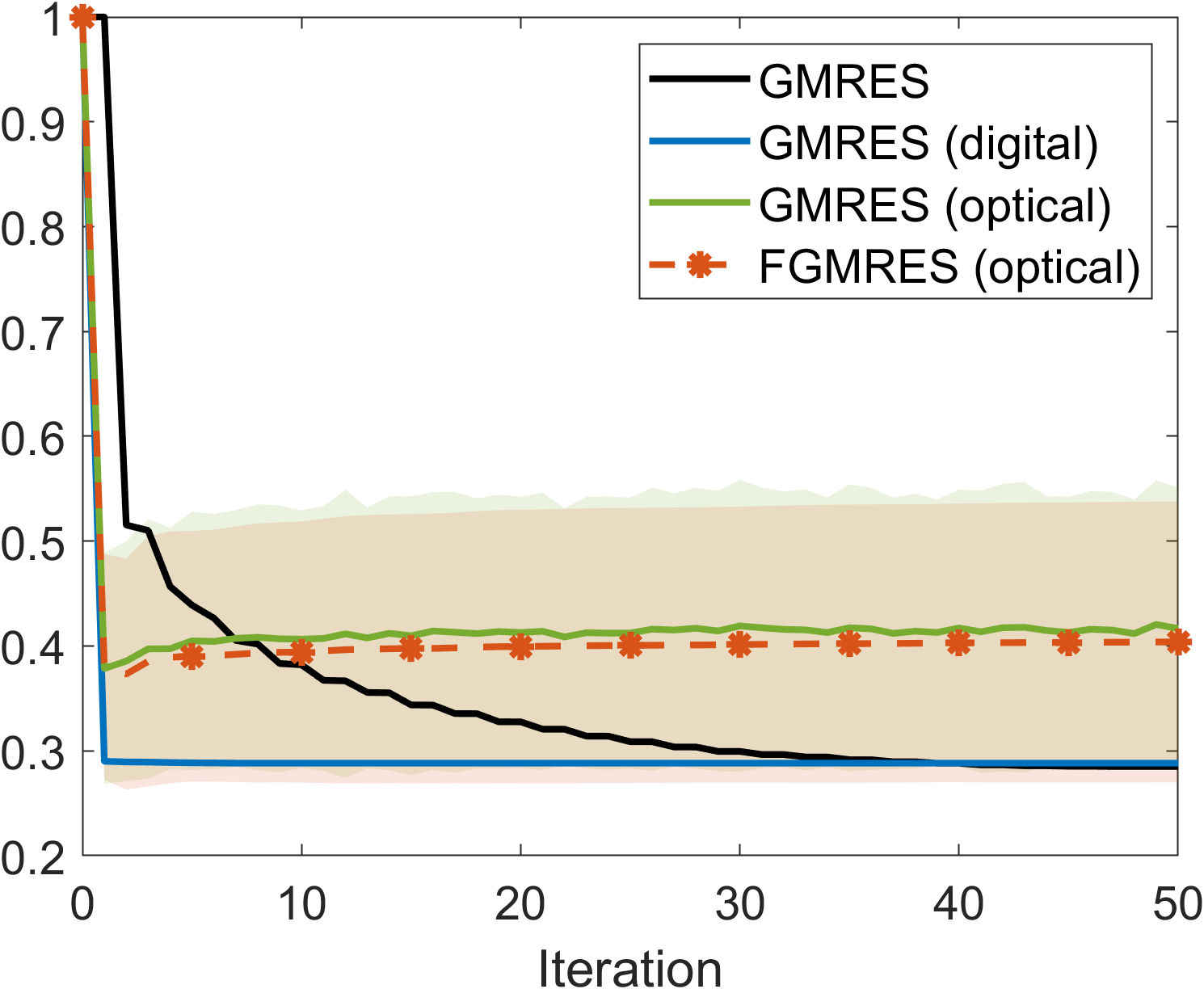}
\end{subfigure}
\caption{RRE across iterations of GMRES/FGMRES with digital and optical Fourier transform-based preconditioning applied to the \emph{Speckle} problem with $5\%$ noise. Optical system qualities of near-ideal, moderate, and degraded are shown from left to right. For the OFT-based methods, results are computed over $50$ independent trials with mean shown by the corresponding curve and the shaded region representing one standard deviation from the mean. Results for GMRES with no preconditioning are provided for reference. 
}
\label{fig:prblurspeck_convergence}
\end{figure}

In Figure~\ref{fig:prblurspeck_convergence}, we present the convergence history for the four methods considered in Table \ref{tab:prblurspeck_results_M1} for the case of $5\%$ noise in $b$. In each plot we present the mean RRE across iterations for all four methods, as well as shaded regions indicating one standard deviation from the mean. The three plots correspond to the near-ideal, moderate, and degraded optical quality levels, from left to right, for the OFT-based methods. Thus, the digital methods considered are identical across all three plots. We again observe comparable performance between GMRES with digital FFT-based preconditioning and FGMRES with OFT-based preconditioning with near-ideal (left plot) and moderate (middle plot) optical system quality. With a degraded optical system quality (right plot), however, we observe larger reconstruction errors and increased variability with the OFT-based methods.

\begin{figure}[h]
\centering
\begin{subfigure}{0.3\textwidth}
    \centering
    \includegraphics[width=\linewidth,height=3.5cm,keepaspectratio]{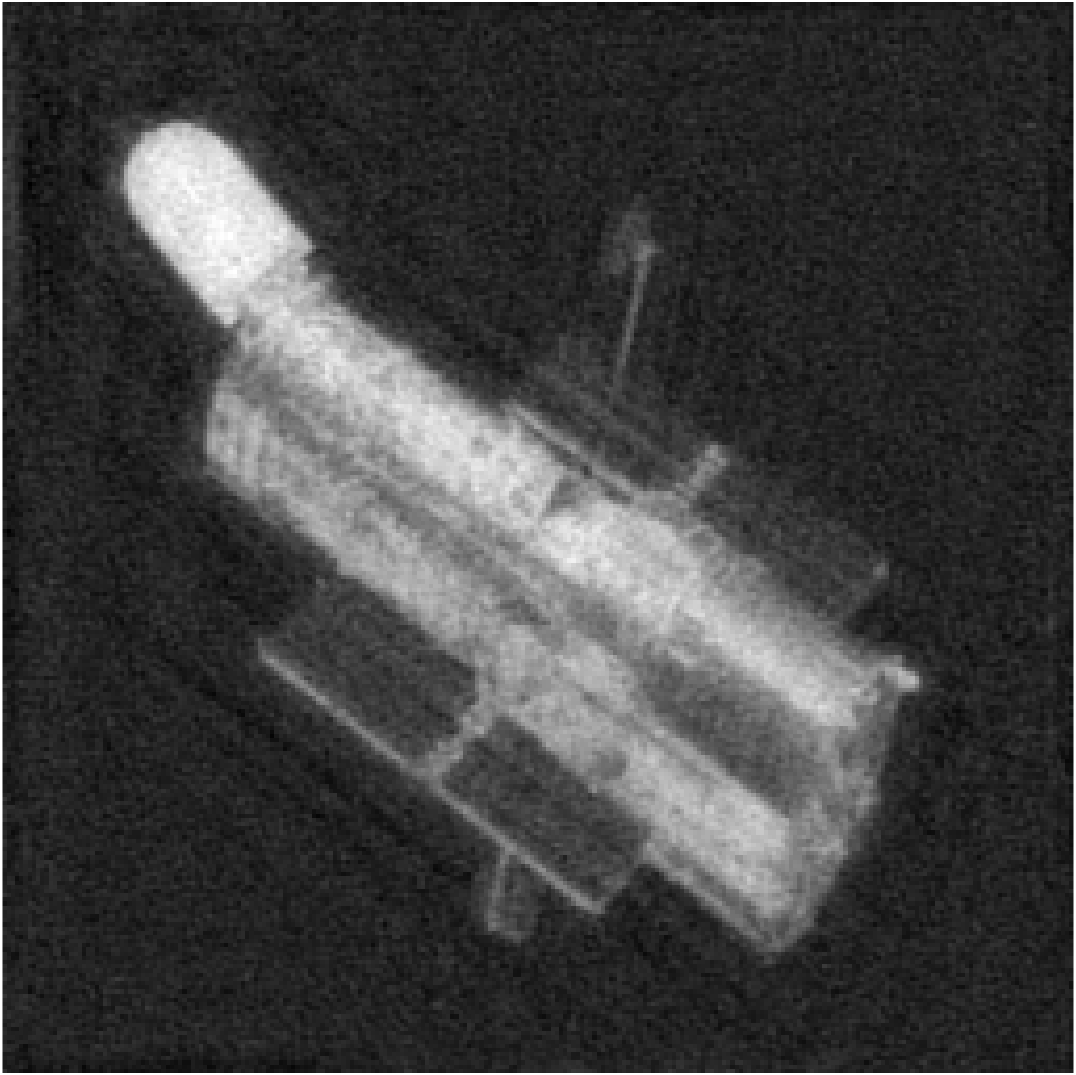}
    \caption{}
    \label{fig:prblurspeck_gmresdig}
\end{subfigure}
\begin{subfigure}{0.3\textwidth}
    \centering
    \includegraphics[width=\linewidth,height=3.5cm,keepaspectratio]{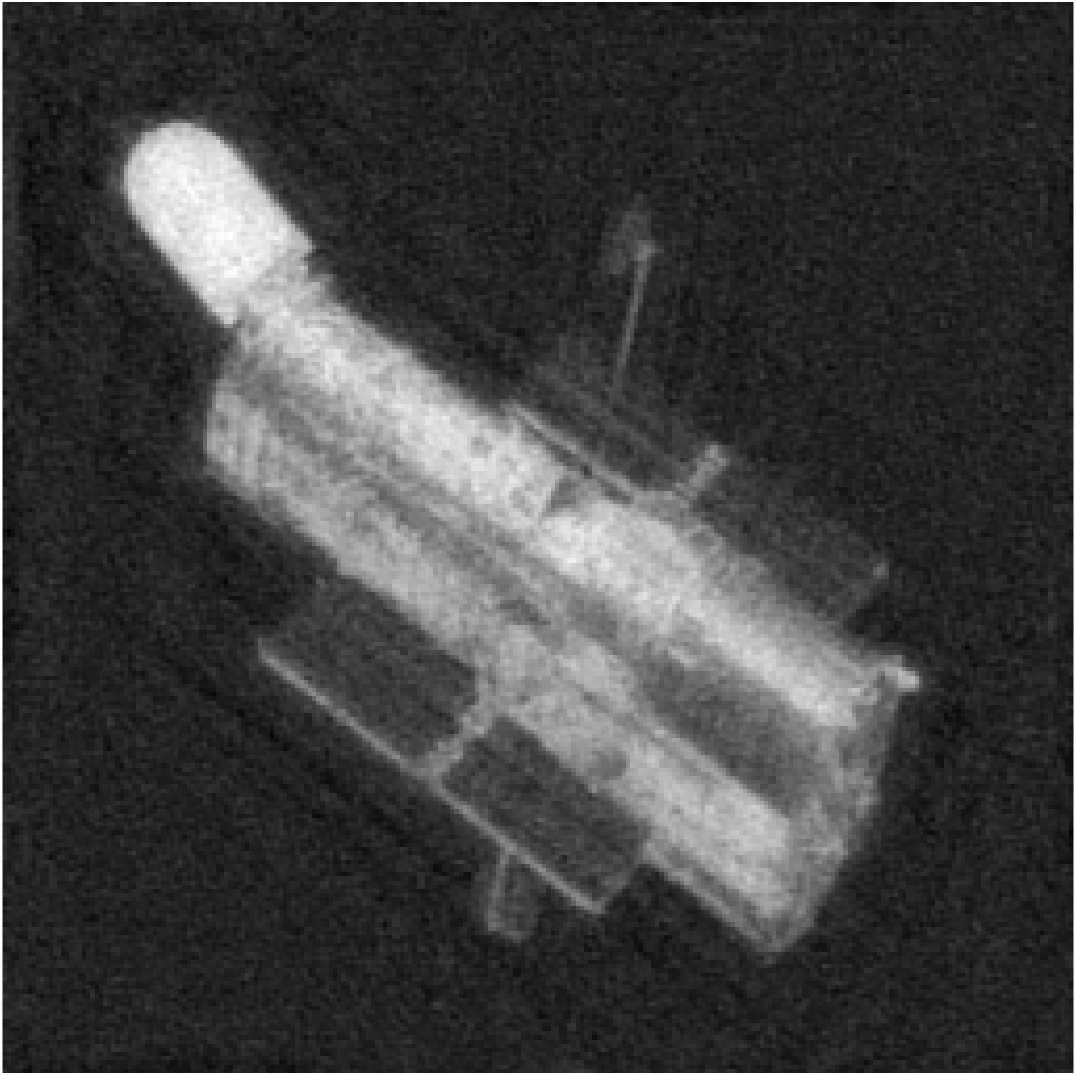}
    \caption{}
    \label{fig:prblurspeck_fgmresopt_nearideal}
\end{subfigure}
\begin{subfigure}{0.3\textwidth}
    \centering
    \includegraphics[width=\linewidth,height=3.5cm,keepaspectratio]{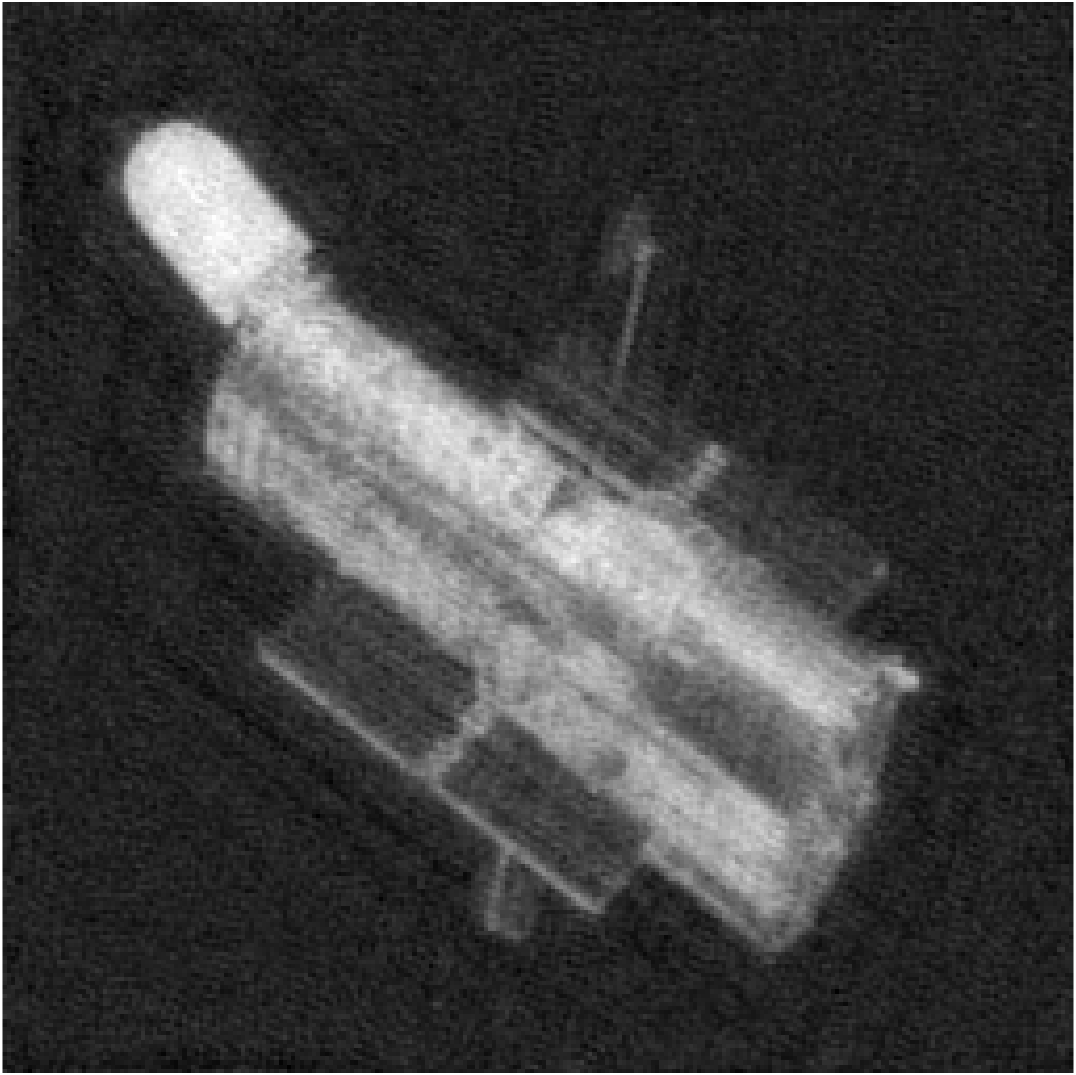}
    \caption{}
    \label{fig:prblurspeck_fgmresopt_moderate}
\end{subfigure}
\vspace{-0.8em}
\caption{Image reconstructions for the \emph{Speckle} test problem with $1\%$ additive Gaussian noise in $b$, produced using the following methods: (a) GMRES with digital FFT-based preconditioning, (b) FGMRES with OFT-based preconditioning (near-ideal optical system quality), and (c) FGMRES with OFT-based preconditioning (moderate optical system quality).}
\label{fig:prblurspeck_recons}
\end{figure}

Lastly, we present a few of the reconstructed images for the case of $1\%$ noise in $b$. The reconstruction produced by GMRES with digital FFT-based preconditioning is shown in Figure~\ref{fig:prblurspeck_gmresdig}. The reconstructions corresponding to the trials with RREs closest to the median RRE of the $50$ trials of FGMRES with OFT-based preconditioning for near-ideal and moderate optical qualities are shown in Figures~\ref{fig:prblurspeck_fgmresopt_nearideal} and \ref{fig:prblurspeck_fgmresopt_moderate}, respectively. We note that the reconstructions seen in these three figures appear indistinguishable to the naked eye, further supporting the comparable reconstruction accuracy between the digitally and optically preconditioned methods.

\subsection{X-Ray Tomography} \label{subsec:prtomo}

\begin{figure}[h]
\centering
\begin{subfigure}{0.3\textwidth}
    \centering
    \includegraphics[height=3.5cm,keepaspectratio]{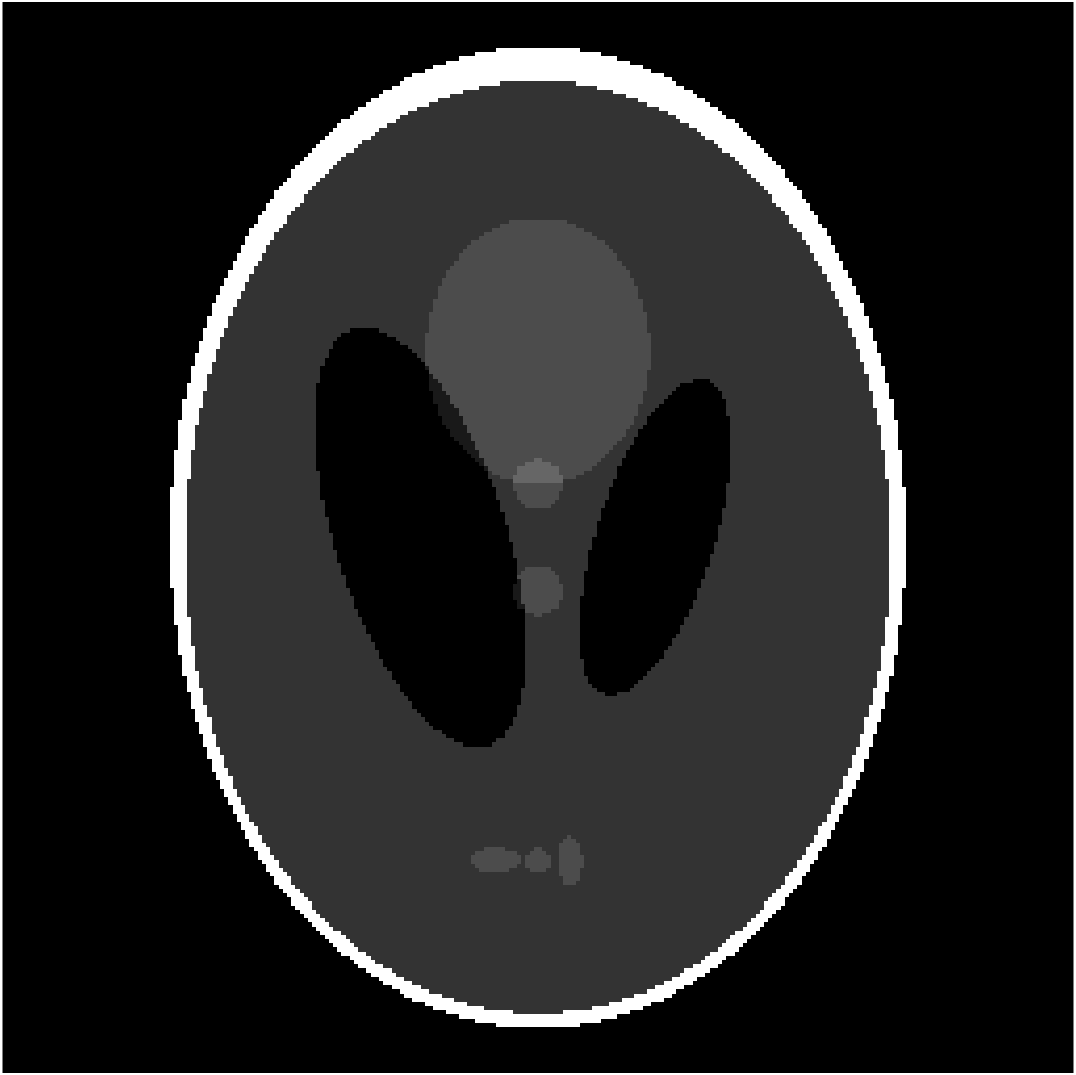}
    \caption{}
    \label{fig:prtomo_xtrue}
\end{subfigure}
\begin{subfigure}{0.3\textwidth}
    \centering
    \includegraphics[height=3.5cm,keepaspectratio]{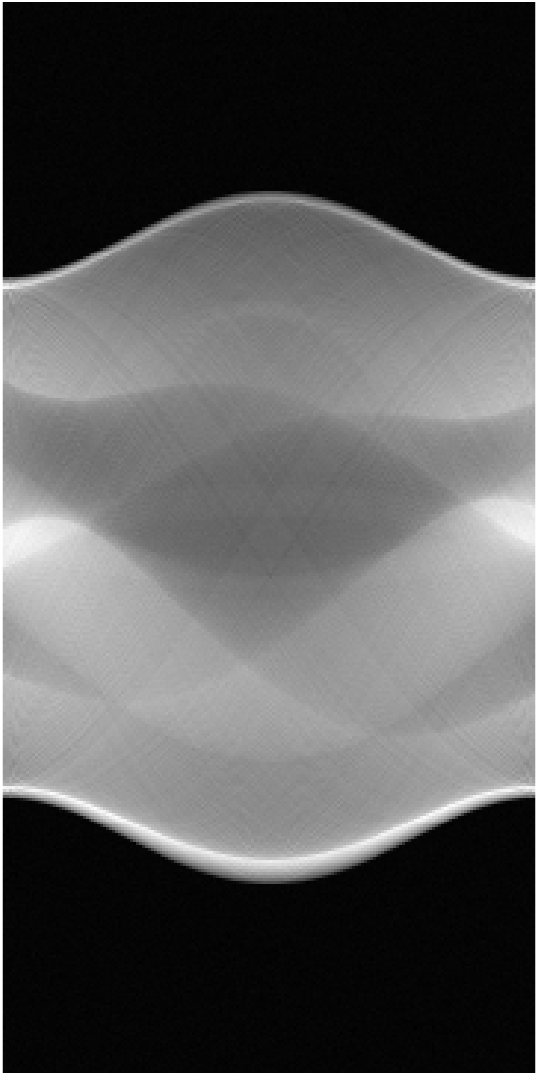}
    \caption{}
    \label{fig:prtomo_bn}
\end{subfigure}
\vspace{-0.8em}
\caption{$2$D X-ray tomography problem \emph{Tomo} with parallel beam geometry: (a) true Shepp-Logan phantom image ($256\times 256$ pixels) and (b) measured sinogram (parallel beam, $180$ projection angles) contaminated by $1\%$ additive Gaussian noise.}
\label{fig:prtomosetup}
\end{figure}

Lastly, we consider the computed tomography test problem from the IR tools package with parallel-beam geometry, where the forward operator $A$ is a discretization of the Radon transform. For our experiments, we consider projection data collected at $180$ equally spaced angles ($0$ to $179$ degrees) and $362$ rays for each source angle. We generate the Shepp-Logan phantom image with $256\times256$ pixels, corresponding to matrix $A$ of size $65160\times65536$. Additive Gaussian noise is included in the right-hand side $b=Ax+e$. In Figures~\ref{fig:prtomo_xtrue} and \ref{fig:prtomo_bn}, the true Shepp-Logan phantom and its corresponding sinogram with $1\%$ additive Gaussian noise is shown. 

Similarly to the \emph{Speckle} problem, we solve the saddle-point formulation shown in (\ref{eq:saddlePoint}) using GMRES and FGMRES with right preconditioning, as described in (\ref{eq:min_saddleProb}). In contrast to the \emph{Speckle} test problem, the tomography forward operator does not readily have a Fourier-diagonalizable approximation of $A$, which is required by the Type 1 preconditioner. Instead, we utilize the Type 2 preconditioner introduced in (\ref{eq:type2precon}), which allows us to exploit an approximation of the normal operator $A^TA$.

For parallel-beam X-ray tomography, the normal operator for the continuous Radon transform $\mathcal{R}$ satisfies
\begin{equation*}
\mathcal{R}^*\mathcal{R} f = 4\pi (-\Delta)^{-1/2}f
\end{equation*}
where $(-\Delta)^{-1/2}$ denotes the inverse square root of the negative Laplacian \cite{Epstein_book}. Motivated by this, we use the following approximation of the discrete normal operator
\begin{equation*}
A^{T}A \approx c(-\Delta)^{-1/2},
\end{equation*}
where $c$ is a scaling constant accounting for discretization. Notably, the fractional Laplacian is diagonalized by the Fourier transform,
\begin{equation*}
(-\Delta)^{-1/2}=F^{-1}\bar{\Lambda} F,
\end{equation*}
where $F$ denotes the discrete Fourier transform and $\bar{\Lambda}$ contains the eigenvalues of the Laplacian operator. This results in the following OFT-compatible (Type 2) preconditioner
\[
K=
\begin{bmatrix}
I & 0\\
A^T & I
\end{bmatrix}
\begin{bmatrix}
I & 0\\
0 & -\left(c(-\Delta)^{-1/2}+\alpha^2I\right)
\end{bmatrix}
\begin{bmatrix}
I & A\\
0 & I
\end{bmatrix}
\]
where the application of $\left(c(-\Delta)^{-1/2}+\alpha^2I\right)^{-1}$ can be implemented using a $4f$ optical system as described in Section \ref{sec:opt_implementation}.

Aside from the preconditioner type and optical implementation, we follow the same experimental setup as the \emph{Speckle} test problem in Section \ref{subsec:prblurspeck}. In particular, we consider the same noise levels for $b$ ($1\%,5\%, \text{and }10\%$), the optical quality levels and parameter ranges described in Table \ref{tab:optparams}, and $50$ trials for the OFT-based methods. We also use the optimal regularization parameter shown in (\ref{eq:regparam}), the discrepancy principle stopping criterion given by (\ref{eq:discrepancy}), and the RRE defined in (\ref{eq:RRE}) for evaluating the reconstruction quality.

In Table \ref{tab:prtomo_results_M3}, we again observe increased reconstruction error as optical system quality decreases, across all OFT-based methods, further demonstrating the importance of sufficient optical system quality in order to maintain the effectiveness of the preconditioner. For all OFT-based cases, the differences between the GMRES and FGMRES results are considerably less pronounced than those seen in \ref{subsec:prblurspeck} for the \emph{Speckle} test problem. Both GMRES, and its flexible variant, with OFT-based preconditioning, have comparable performance to the digital FFT-based version, in terms of both reconstruction error and iteration count for the near-ideal and moderate optical quality ranges across all noise levels considered in $b$. We note, however, that for the $1\%$ noise level case, the preconditioned methods are unable to achieve the reconstruction error seen without preconditioning.

\begin{figure}[h]
\centering
\begin{subfigure}[h]{0.46\columnwidth}
    \centering
    \includegraphics[width=\linewidth]{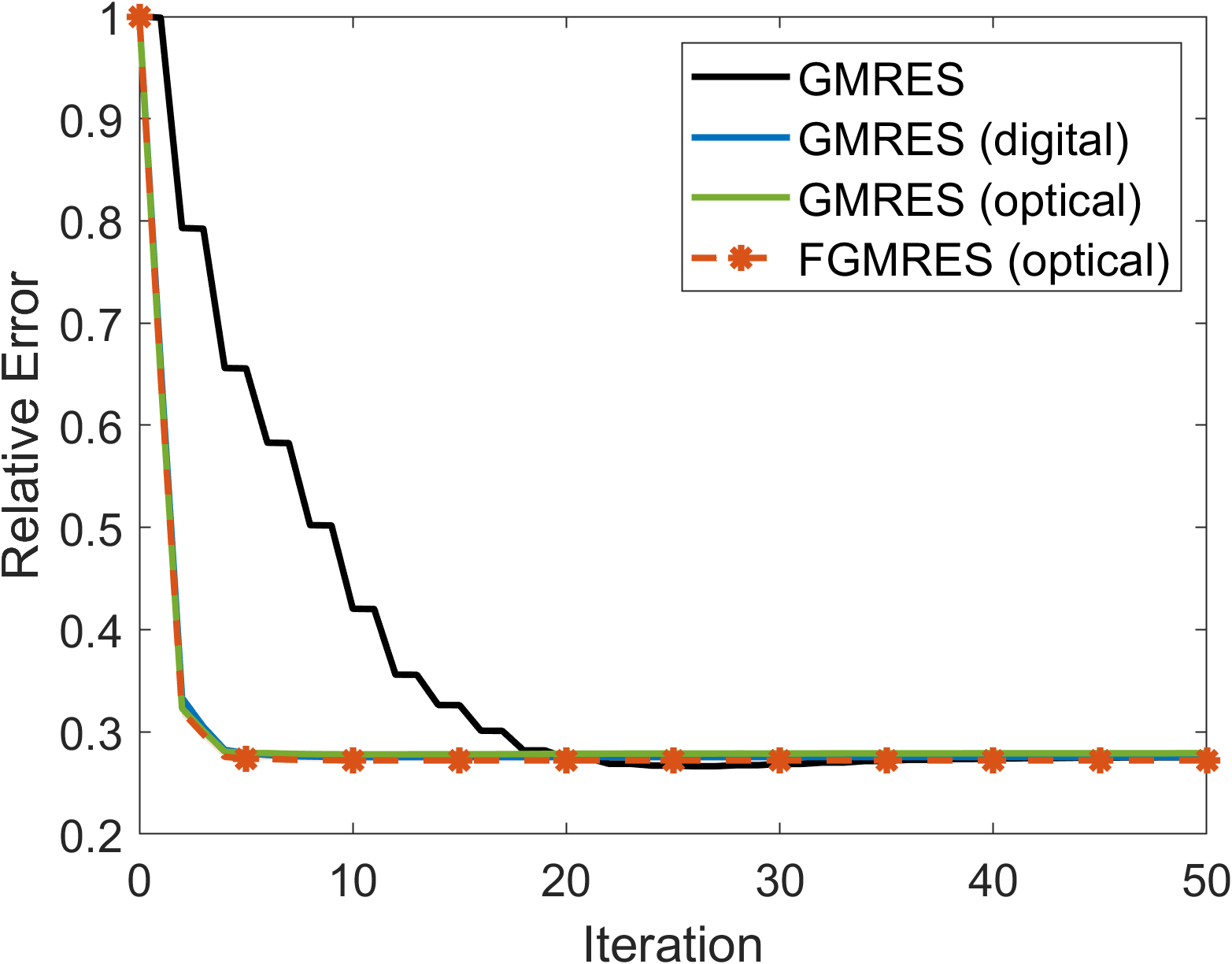}
\end{subfigure}
\hspace{.5cm}
\begin{subfigure}[h]{0.46\columnwidth}
    \centering
    \includegraphics[width=\linewidth]{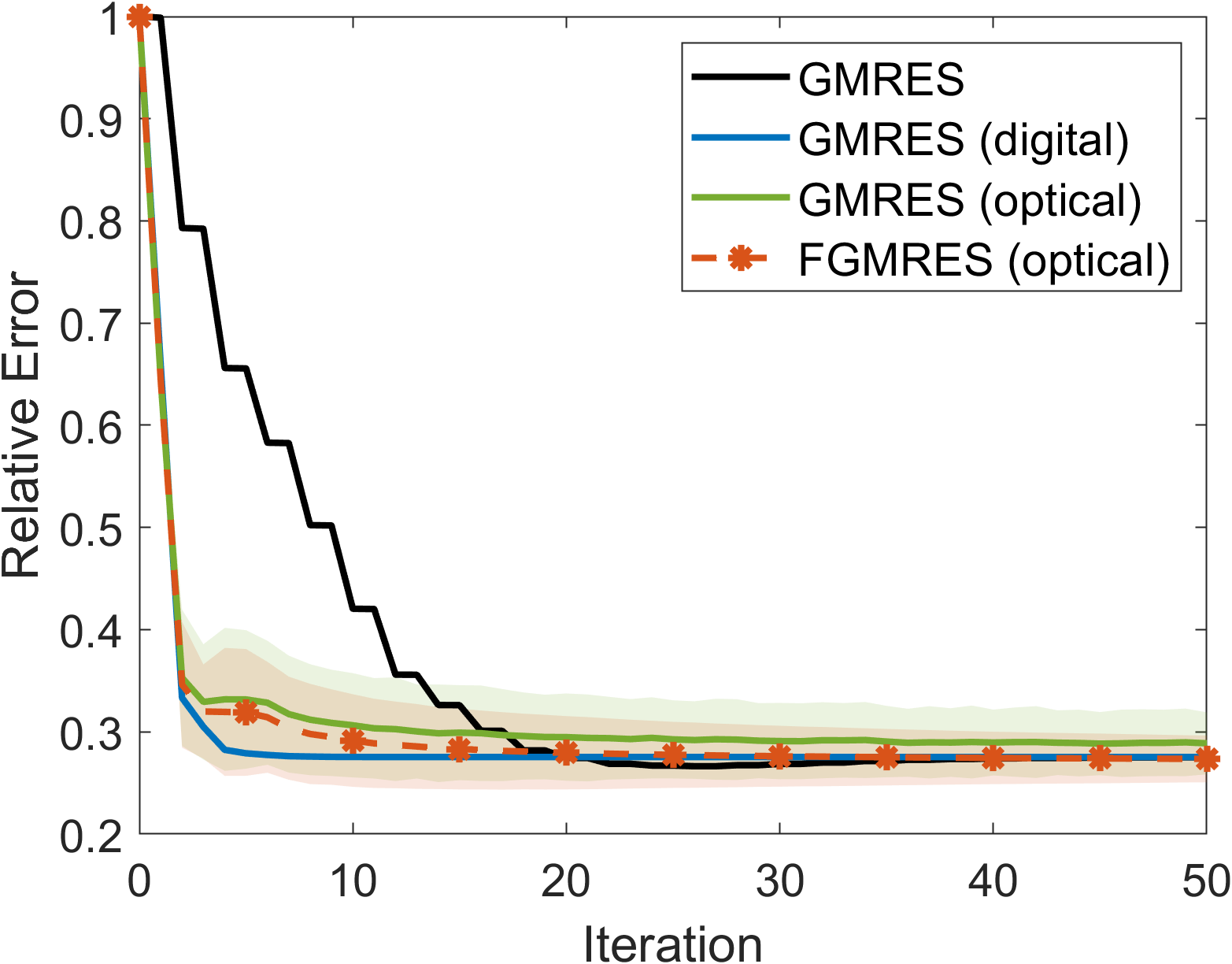}
\end{subfigure}
\caption{Relative reconstruction error across iterations of GMRES/FGMRES with digital and optical Fourier transform-based preconditioning applied to the \emph{Tomo} problem with $5\%$ noise. Optical system qualities of near-ideal and moderate are shown in the left and right plots, respectively. For the OFT-based methods, results are computed over $50$ independent trials with mean shown by the corresponding curve and the shaded region representing one standard deviation from the mean. Results for GMRES with no preconditioning are provided for reference.}
\label{fig:prblurtomo_convergence}
\end{figure}

In Figure~\ref{fig:prblurtomo_convergence}, we present the convergence history for the four methods considered in Table \ref{tab:prtomo_results_M3} for the case of $5\%$ noise in $b$. In each plot we present the mean RRE across iterations for all four methods, as well as shaded regions indicating one standard deviation from the mean. The two plots correspond to the near-ideal (left) and moderate (right) optical quality levels for the OFT-based methods. The digitally and optically preconditioned methods have comparable reconstruction performance. However, it can be seen that with moderate optical system quality (right plot) there may be slower convergence than in the near-ideal and digital cases. 

\begin{figure}[h]
\centering
\begin{subfigure}{0.3\textwidth}
    \centering
    \includegraphics[height=3.5cm,keepaspectratio]{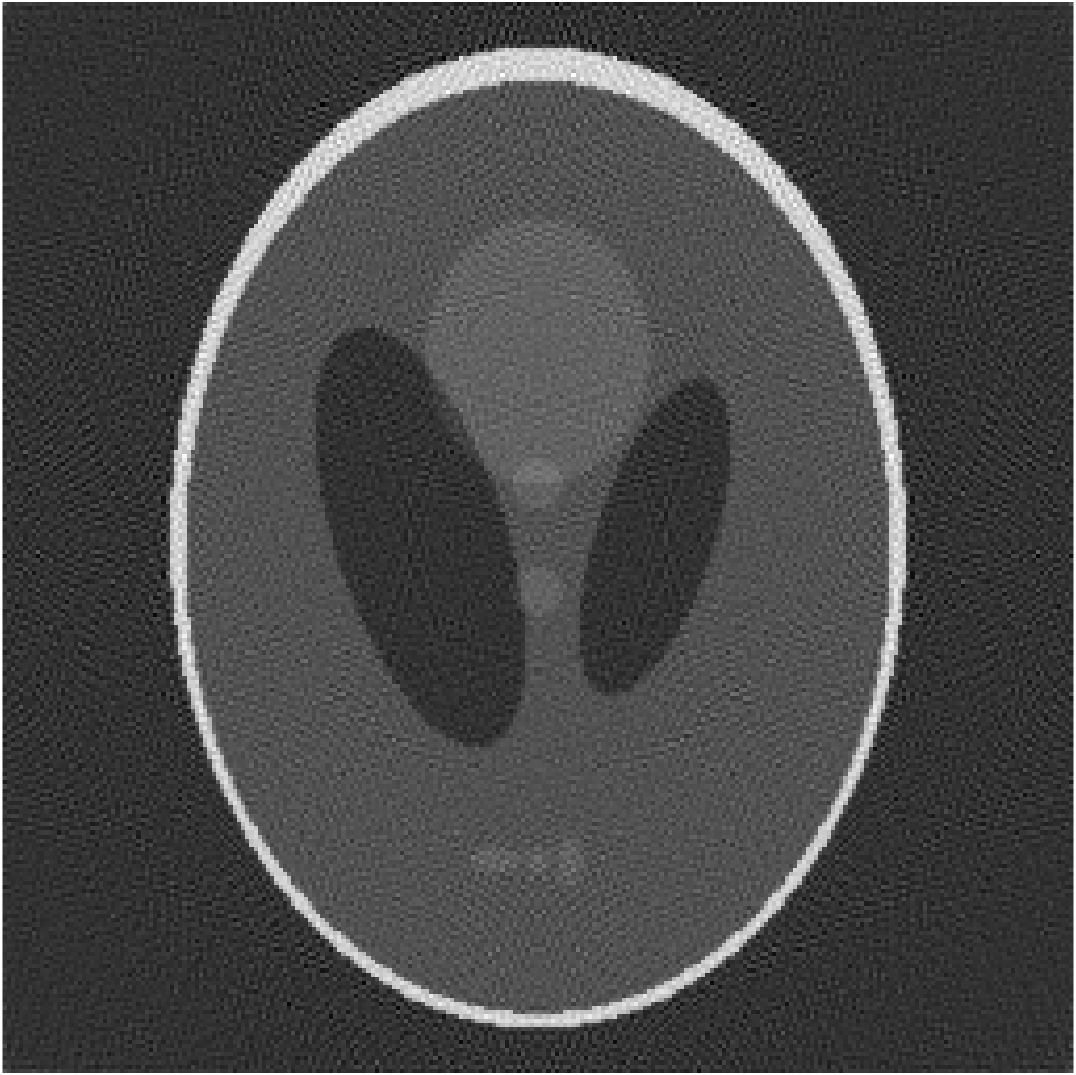}
    \caption{}
    \label{fig:prtomo_recondig}
\end{subfigure}
\begin{subfigure}{0.3\textwidth}
    \centering
    \includegraphics[height=3.5cm,keepaspectratio]{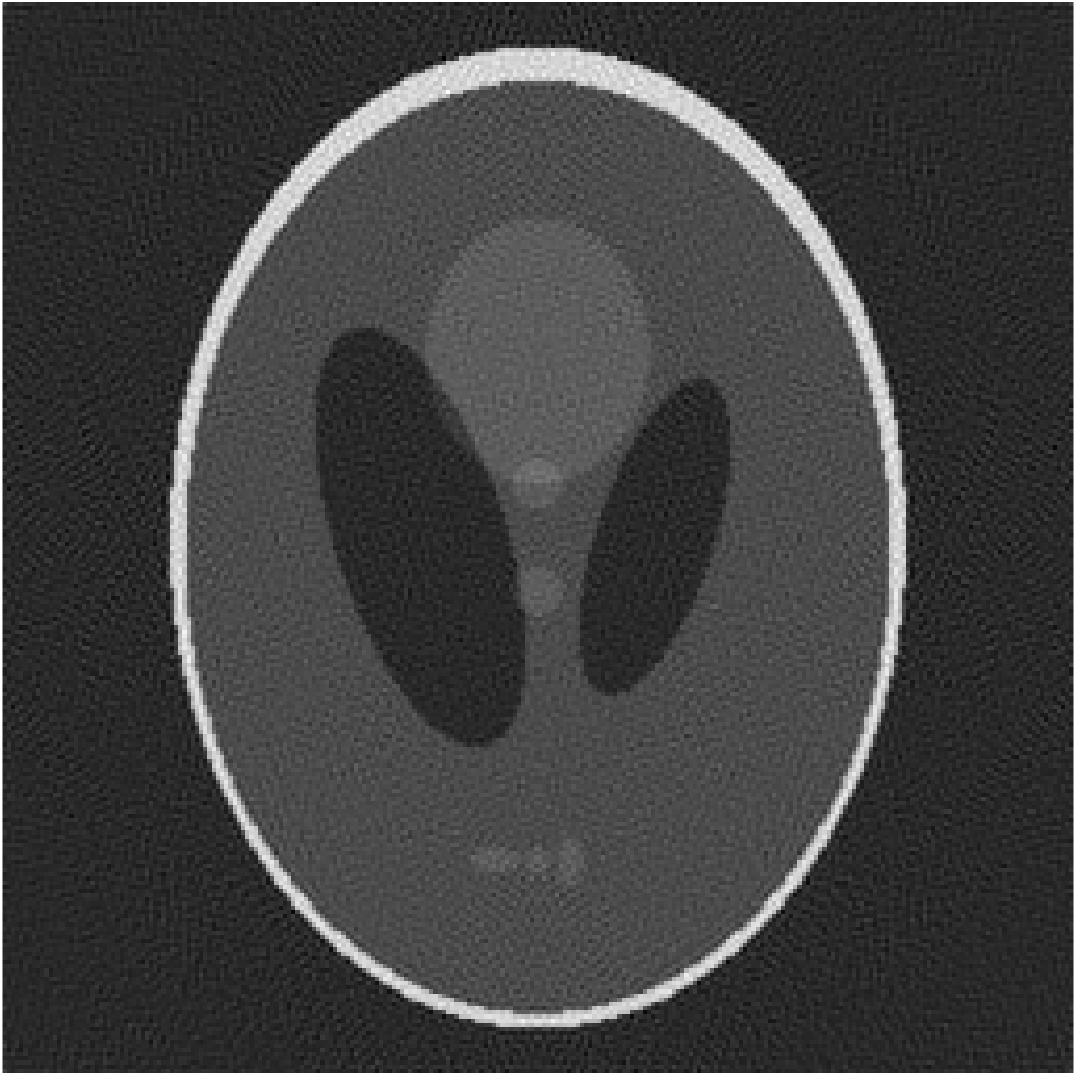}
    \caption{}
    \label{fig:prtomo_recon_nearideal}
\end{subfigure}
\begin{subfigure}{0.3\textwidth}
    \centering
    \includegraphics[height=3.5cm,keepaspectratio]{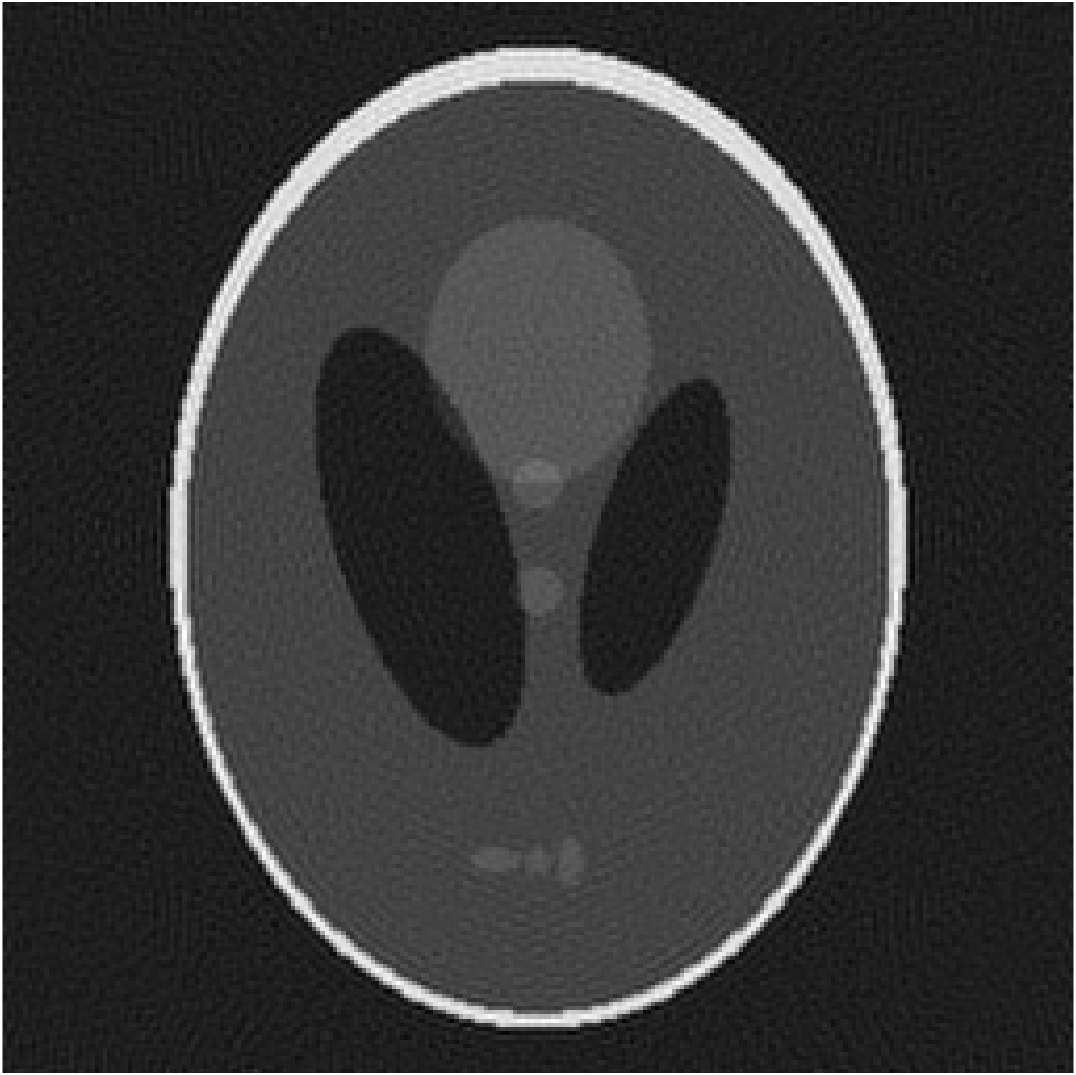}
    \caption{}
    \label{fig:prtomo_recon_moderate}
\end{subfigure}
\vspace{-0.8em}
\caption{Reconstructed images for the $2$D X-ray tomography problem \emph{Tomo} with parallel beam geometry and $1\%$ additive Gaussian noise in $b$, produced using the following methods: (a) GMRES with digital FFT-based preconditioning, (b) FGMRES with OFT-based preconditioning (near-ideal optical system quality), and (c) FGMRES with OFT-based preconditioning (moderate optical system quality).}
\label{fig:prtomorecons}
\end{figure}

We conclude this section by presenting a few of the reconstructed images for the case of $1\%$ noise in $b$. The reconstruction produced by GMRES with digital FFT-based preconditioning is shown in Figure~\ref{fig:prtomo_recondig}. The reconstructions corresponding to the trials with RREs closest to the median RRE of the $50$ trials of FGMRES with OFT-based preconditioning for near-ideal and moderate optical qualities are shown in Figures~\ref{fig:prtomo_recon_nearideal} and \ref{fig:prtomo_recon_moderate}, respectively. All three reconstructions appear indistinguishable, consistent with the comparable reconstruction accuracy between the digitally and optically preconditioned methods.

\begin{table}[h]
\centering \footnotesize
\renewcommand{\arraystretch}{1.3}
\setlength{\tabcolsep}{4pt}
\caption{RRE results for the \emph{Tomo} test problem with $1\%$, $5\%$, and $10\%$ noise. GMRES with no preconditioning, GMRES with digital FFT-based preconditioning, GMRES with OFT-based preconditioning, and FGMRES with OFT-based preconditioning are considered. For the OFT-based methods, the optical system quality is varied among near-ideal, moderate, and degraded. Results for the OFT-based methods are computed over $50$ independent trials. The mean and standard deviation of the RRE and the corresponding iteration count is reported. Also reported for the OFT-based methods is the minimum and maximum RRE with their corresponding iteration counts. The RRE values presented correspond to the $p^{\text{th}}$ iterative solution when the discrepancy principle is satisfied, or alternatively, when the specified maximum number of iterations is achieved ($p=100$).} 
\begin{tabular}{c | l | c | l | l | l}
\multirow{2}{*}{\makecell{\textbf{Noise}\\\textbf{in $\boldsymbol{b}$}}} & \multirow{2}{*}{\textbf{Method}} & \multirow{2}{*}{\makecell{\textbf{Optical}\\\textbf{Quality}}} & \multicolumn{3}{c}{\textbf{RRE}} \\
\cline{4-6}
& & & \textbf{$\boldsymbol{\mu \pm \sigma}$ (Iter. $\boldsymbol{\mu \pm \sigma}$)} &
\textbf{Min (Iter.)} &
\textbf{Max (Iter.)} \\
\Xhline{1.5pt} 

\multirow{8}{*}{$1\%$} & GMRES & --- & 0.1680 (32) & \hspace{.6cm}--- & \hspace{.6cm}--- \\
\cline{2-6}
& GMRES (dig.) & --- & 0.2811 (6) & \hspace{.6cm}--- & \hspace{.6cm}--- \\
\cline{2-6}
& \multirow{3}{*}{GMRES (opt.)} & Near-ideal & $0.2055 \pm 0.0036$ ($100 \pm 0$) & 0.1986 (100) & 0.2170 (100) \\
& & Moderate & $0.2558 \pm 0.0690$ ($100 \pm 0$) & 0.2130 (100) & 0.6590 (100) \\
& & Degraded & $1.0628 \pm 1.3688$ ($100 \pm 0$) & 0.3345 (100) & 7.7480 (100) \\
\cline{2-6}
& \multirow{3}{*}{FGMRES (opt.)} & Near-ideal & $0.2069 \pm 0.0245$ ($28 \pm 39$) & 0.1726 (14) & 0.2408 (6) \\
& & Moderate & $0.1954 \pm 0.0541$ ($100 \pm 0$) & 0.1746 (100) & 0.5563 (100) \\
& & Degraded & $0.8786 \pm 1.1298$ ($100 \pm 0$) & 0.2146 (100) & 5.0176 (100) \\
\hline \hline 

\multirow{8}{*}{$5\%$} & GMRES & --- & 0.2745 (20) & \hspace{.6cm}--- & \hspace{.6cm}--- \\
\cline{2-6}
& GMRES (dig.) & --- & 0.3049 (3) & \hspace{.6cm}--- & \hspace{.6cm}--- \\
\cline{2-6}
& \multirow{3}{*}{GMRES (opt.)} & Near-ideal & $0.2807 \pm 0.0024$ ($4 \pm 0$) & 0.2768 (4) & 0.2849 (4) \\
& & Moderate & $0.3172 \pm 0.0320$ ($8 \pm 19$) & 0.2800 (4) & 0.4275 (100) \\
& & Degraded & $0.7169 \pm 0.8295$ ($92 \pm 24$) & 0.3288 (9) & 4.4355 (100) \\
\cline{2-6}
& \multirow{3}{*}{FGMRES (opt.)} & Near-ideal & $0.2759 \pm 0.0035$ ($4 \pm 0$) & 0.2718 (4) & 0.2957 (3) \\
& & Moderate & $0.3082 \pm 0.0321$ ($6 \pm 12$) & 0.2705 (4) & 0.3842 (86) \\
& & Degraded & $0.6377 \pm 0.6765$ ($69 \pm 39$) & 0.3199 (58) & 3.3747 (100) \\
\hline \hline 

\multirow{8}{*}{$10\%$} & GMRES & --- & 0.3728 (14) & \hspace{.6cm}--- & \hspace{.6cm}--- \\
\cline{2-6}
& GMRES (dig.) & --- & 0.3731 (3) & \hspace{.6cm}--- & \hspace{.6cm}--- \\
\cline{2-6}
& \multirow{3}{*}{GMRES (opt.)} & Near-ideal & $0.3652 \pm 0.0020$ ($3 \pm 0$) & 0.3617 (3) & 0.3701 (3) \\
& & Moderate & $0.3796 \pm 0.0260$ ($3 \pm 0$) & 0.3609 (3) & 0.5107 (3) \\
& & Degraded & $0.6007 \pm 0.4945$ ($49 \pm 45$) & 0.3853 (100) & 2.6776 (100) \\
\cline{2-6}
& \multirow{3}{*}{FGMRES (opt.)} & Near-ideal & $0.3642 \pm 0.0020$ ($3 \pm 0$) & 0.3608 (3) & 0.3690 (3) \\
& & Moderate & $0.3767 \pm 0.0222$ ($3 \pm 0$) & 0.3584 (3) & 0.4818 (3) \\
& & Degraded & $0.5738 \pm 0.4394$ ($47 \pm 46$) & 0.3806 (100) & 2.4743 (100) \\
\hline 
\end{tabular}
\label{tab:prtomo_results_M3}
\end{table}

\section{Conclusions} \label{sec:conclusions} 

In this work we proposed the use of OFT-based preconditioners applied to certain GMRES-type iterative methods for the solution of discrete inverse problems. Due to the inherent errors present in practical OFTs, we investigated the feasibility of this hybrid digital-optical approach. We simulated the OFT digitally by incorporating sources of error expected in practical optical systems, including finite aperture, aberrations, misalignment, noise, and quantization. Our numerical results suggest that, given sufficient optical system quality, iterative methods using OFT-based preconditioners can remain effective, despite the errors being introduced into the system by the OFT. Specifically, the simulated hybrid digital-optical methods achieved reconstruction quality and convergence similar to those of their fully digital counterparts. Our results demonstrate the potential effectiveness of OFT-based preconditioners for certain discrete inverse problems, thereby warranting their continued study in physically implemented hybrid digital-optical systems. Additional future work should investigate the computational complexity and speedup of the proposed hybrid digital-optical preconditioned strategy as well as perform dedicated experimental validation studies.


\section*{Acknowledgments}
This work was partially funded by the U.S. National Science Foundation, under grant DMS-2038118. Any opinions, finding, and conclusions or recommendations expressed in this material are those of the author(s) and do not necessarily reflect the views of the National Science Foundation. This research used resources of the Argonne Leadership Computing Facility at Argonne National Laboratory, which is supported by the Office of Science of the U.S. Department of Energy, Office of Science, under contract number DE-AC02-06CH11357.

\bibliographystyle{siamplain}
\bibliography{references}
\end{document}